\documentclass[12pt]{amsart}
\usepackage{amsmath,amsthm,mathrsfs,amssymb,wasysym}
\usepackage{geometry,hyperref}

\newtheorem{theorem}{Theorem}[section]

\newtheorem{lemma}[theorem]{Lemma}
\newtheorem{proposition}[theorem]{Proposition}

\theoremstyle{definition}

\begin{document}
\title[Symmetric square $L$-functions]{A zero-density estimate for the symmetric square $L$-functions in the weight aspect}
\author[S. Yin]{Songchen Yin}
\address{Data Science Institute and School of Mathematics\\ Shandong University\\
 Jinan, Shandong 250100, China}
\email{scyin@mail.sdu.edu.cn}
\date{}
\maketitle

\begin{abstract}
In this paper, we prove a zero-density estimate for the symmetric square $L$-functions associated with Hecke eigencuspforms in the weight aspect. The proofs are based on the mollified second moment of these $L$-functions.
\end{abstract}

\frenchspacing

\textwidth=13.5cm
\textheight=23cm
\parindent=16pt
\oddsidemargin=-0.5cm
\evensidemargin=-0.5cm
\topmargin=-0.5cm
\linespread{1.2}

\section{Introduction}

The zero-density estimate is a classic problem in the analytic theory of $L$-functions. It is one of the ways to avoid the Riemann Hypothesis in applications. Roughly speaking, the aim of the zero-density estimate is to show that ``most'' zeros lie not too far from the critical line $\Re s=1/2$. Generally, for an $L$-function $L(s,f)$ absolutely converge in $\Re s>1$, if it satisfies the functional equation and admits analytic continuation to $\mathbb{C}$, we define
\[
N(\sigma,T,f)=\#\{\rho=\beta+i\gamma:L(\rho,f)=0,\beta>\sigma,0<\gamma<T\},
\]
and write $N(\sigma,T)=N(\sigma,T,\zeta)$.

A landmark result in this direction is due to Selberg \cite{selberg1946contribution}, who proved
\[
N(\sigma,T)=O\left(T^{1-\frac{1}{4}(\sigma-\frac{1}{2})}\log T\right).
\]
This estimate has two distinctive features. First, the power of $\log T$
matches the order of magnitude of the $\log T$ term in the asymptotic formula
for $N(T)=N(0,T)$,
\[
N(T)=\frac{T}{2\pi}\log\frac{T}{2\pi}-\frac{T}{2\pi}+O(\log T).
\]
Second, it remains non-trivial even when $\sigma-\frac{1}{2}$ is as small as
$c/\log T$ for some large constant $c$. This result is also a significant step in proving the Selberg's central limit theorem, i.e. the real and imaginary part of the Riemann zeta function behave as a Gaussian random variable asymptotically after normalization. Selberg \cite{selberg1946contributions} also established an
analogous zero-density estimate for the family of Dirichlet $L$-functions
$L(s,\chi)$ with $\chi$ ranging over characters modulo $q$.

For $L$-functions attached to modular forms, Luo \cite{luo1995zeros} proved
a Selberg-type zero-density estimate in the $t$-aspect for $L$-functions
associated with Hecke eigencuspforms for $\mathrm{SL}_2(\mathbb{Z})$:
\[
	N(\sigma,T,f)=O_f\left(T^{1-\frac{1}{72}(\sigma-\frac{1}{2})}\log T\right).
\]
Beckwith--Liu--Thorner--Zaharescu  \cite{beckwith2023zero} improved the
constant $\frac{1}{72}$ to $\frac{25}{128}$ and extended the result to
Hecke--Maass forms. 

Hough \cite{hough2012zero} proved a corresponding result in
the weight aspect for forms of large weight $k$:
\[
	\frac{1}{|H_k|}\sum_{f\in H_k}N(\sigma,T,f)=O\left(Tk^{-\theta\left(\sigma-\frac{1}{2}\right)}\log k\right),
\]
uniformly for $\frac{1}{2}+\frac{2}{\log k}<\sigma<1$ and $\frac{10}{\log k}<T<k^{\delta}$, with $\theta$ and $\delta$ sufficiently small. Kowalski--Michel \cite{kowalski1999analytic} established an analogous estimate in the level aspect for forms of weight $2$ and large level $q$. Recently, Sun--Wang--Yu \cite{sun2025selberg} established a corresponding result for the family of
twisted $\mathrm{GL}(2)$ $L$-functions, with $\chi$ ranging over primitive characters modulo $q$, and Sun--Wang \cite{sun2024zero} proved a weighted zero-density estimate for $\mathrm{GL}(3)$ $L$-functions in the spectral aspect.

The aim of this article is to extend Hough's result to a family of symmetric
square $L$-functions, working with
$\sigma-\frac{1}{2}>\frac{1}{\log K}$ and $|t|<K^{\delta}$ for some small
$\delta>0$. We begin by recalling the relevant definitions. Suppose $f$ is a
Hecke eigencuspform for $\mathrm{SL}_2(\mathbb{Z})$ of weight $k$, with
Fourier expansion
\[
f(z)=\sum_{n=1}^\infty\lambda_f(n)n^{\frac{k-1}{2}}e(nz).
\]
here we normalize the coefficient such that $\lambda_f(1)=1$, in our normalization $\lambda_f$ satisfies the Deligne's bound \cite{deligne1974conjecture} $|\lambda_f(n)|\leq d(n)$.
If $\Re s>1$, we define the symmetric square $L$-function as an absolutely convergent Dirichlet series:
\[
L(s,\mathrm{sym}^2f)=\zeta(2s)\sum_{n=1}^{\infty}\frac{\lambda_f(n^2)}{n^s},
\]
which has the Euler product
\[L(s,\mathrm{sym}^2f)=\prod_p \left(1-\frac{\lambda_f(p^2)}{p^s}+\frac{\lambda_f(p^2)}{p^{2s}}-\frac{1}{p^{3s}}\right)^{-1}. \]
Therefore we can get the Euler product of the inverse of $L(s,\mathrm{sym}^2f)$ for $\Re(s)>1$:\\
\[
L(s,\mathrm{sym}^2f)^{-1}=\prod_p \left(1-\frac{\lambda_f(p^2)}{p^s}+\frac{\lambda_f(p^2)}{p^{2s}}-\frac{1}{p^{3s}}\right)=\sum_{n=1}^{\infty}\frac{b_f(n)}{n^s}, 
\]
where $b_f(n)$ is a multiplicative function satisfying
\[
b_f(p^l) = 
\begin{cases}
    1 & \text{if } l=0, \\
    -\lambda_f(p^2) & \text{if } l=1, \\
    \lambda_f(p^2) & \text{if } l=2, \\
    -1 & \text{if } l=3, \\
    0 & \text{if } l\geq 4.
\end{cases}
\]

Shimura \cite{shimura1975holomorphy} proved that
$L(s,\mathrm{sym}^2f)$ admits analytic continuation to an entire function on $\mathbb{C}$ and satisfies the functional equation
\[
L_\infty(s,\mathrm{sym}^2f)L(s,\mathrm{sym}^2f)=L_\infty(1-s,\mathrm{sym}^2f)L(1-s,\mathrm{sym}^2f),
\]
where
\[
L_\infty(s,\mathrm{sym}^2f)=\pi^{-\frac{3}{2}s+\frac{1}{2}}2^{2-s-k}\Gamma\left(\frac{s+1}{2}\right)\Gamma(s+k-1).
\]

In the following text, we will extensively deal with the sum over $f\in H_k$, then we need to introduce the harmonic sum $\sum_{f\in H_k}^h$. For any complex numbers $\alpha_f$, let
\[
\sideset{}{^h}\sum_{f\in H_k}\alpha_f=\sum_{f\in H_k}w_f\alpha_f,
\]
be the harmonic average of $\alpha_f$, where
\[
w_f=\frac{12\zeta(2)}{(k-1)L(1,\mathrm{sym}^2f)},
\]
is the weight arising from the Petersson trace formula. We are now in a
position to state our main theorem.

\begin{theorem}\label{T:1}
 Let $N(\sigma,t_1,t_2,\mathrm{sym}^2 f)$ denotes the number
\[
	\#\left\{\rho=\beta+i\gamma:L(\rho,\mathrm{sym}^2 f)=0,\beta>\sigma,t_1<\gamma<t_2\right\},
\]
then for $0<t_1<t_2< K^{\frac{1}{5}}$, $t_2-t_1>1/\log K$, $1/2+1/\log K<\sigma<1$ and $\psi$ is a smooth function with support in $(0,\infty)$, 
\[
	\sum_{k\equiv 0 \bmod 2}\psi\left(\frac{k-1}{K}\right)\sideset{}{^h}\sum_{f\in H_k}N(\sigma,t_1,t_2,\mathrm{sym}^2 f)=O\left((t_2-t_1)K^{1-\frac{1}{21}(\sigma-\frac{1}{2})}\log K\right).
\]
\end{theorem}

Our proof returns to Selberg's original argument \cite{selberg1946contributions} for Dirichlet $L$-functions, which reduces the zero-density estimate to an upper bound for a mollified second moment. Accordingly, a key step is to
estimate the mollified second moment of $L(s,\mathrm{sym}^2f)$. Hough \cite{hough2012zero} obtained such an upper bound for $L(s,f)$. Following his approach, we adapt the mollification method to the symmetric square setting. The mollification is preceded by an evaluation of the twisted second moment: Khan
\cite{khan2010non} established an asymptotic formula for the twisted second moment of $L(1/2,\mathrm{sym}^2f)$, and we adapt his method to treat $L(\sigma+it,\mathrm{sym}^2f)$ for general $\sigma+it$.

The first ingredient is the following asymptotic formula for the twisted
second moment of symmetric square $L$-functions in the weight aspect.

\begin{proposition}[Twisted second moment]\label{P:1} Assume $0<|t|\leq K^{\delta}$ with $0<\delta<1/4$, $r<K^{1-\delta-\epsilon}$, $1/2<\sigma<3/2-\epsilon$ and $\psi$ is a smooth compactly supported function. Take $d=d_1d_2^2$ with $d_1$ squarefree. We have
\[\sum_{k\equiv 0 \bmod 2}\psi\left(\frac{k-1}{K}\right)\sideset{}{^h}\sum_{f\in H_k} |L(\sigma+it,\mathrm{sym}^2f)|^2\lambda_f(r^2)= S_++S_-+S_{OD}+R(r). \\
    \]
Where 
\[
S_+= \sum_{k\equiv 0 \bmod 2}\psi\left(\frac{k-1}{K}\right)\sum_{d|r^2}\frac{1}{(d_1r)^{\sigma}(r/d)^{it}}\zeta(2\sigma+2it)\zeta(2\sigma-2it)\zeta(2\sigma),
\]
\begin{align*}
S_-&=\sum_{k\equiv 0 \bmod 2}\psi\left(\frac{k-1}{K}\right)\sum_{d|r^2}\frac{1}{(d_1r)^{1-\sigma}(r/d)^{it}}\left(\frac{k^2}{4\pi^3}\right)^{1-2\sigma}G_-(0,\sigma,t)\\
&\quad\times\zeta(1-2\sigma+2it)\zeta(1-2\sigma-2it)\zeta(1-2\sigma),
\end{align*}
\begin{align*}
S_{OD}&=\Re\frac{K^{2-2\sigma+2it}}{(4\pi^3)^{1/2-\sigma+ it}}\frac{\Gamma\left(\frac{2-\sigma+it}{2}\right)}{\Gamma\left(\frac{\sigma+1-it}{2}\right)}\tilde{\psi}(-2\sigma+2it)\sum_{d|r^2}\frac{(r/d)^{\sigma-1/2}}{(d_1r)^{1/2+it}}\\
&\quad\times\zeta(2\sigma+2it)\zeta(1+2it)\zeta(2-2\sigma+2it).
\end{align*}
Here $\tilde{\psi}$ is the Mellin transform of $\psi$, $G(s)=e^{s^2}$and
\[
G_-(0,\sigma,t)=\frac{\Gamma((2-\sigma+it)/2)\Gamma((2-\sigma-it)/2)}{\Gamma((1+\sigma+it)/2)\Gamma((1+\sigma-it)/2)},
\]
\[
R(r)\ll K^{\epsilon}r^{2\sigma-1}+K^{-\sigma+\epsilon}r^{1-\sigma}|t|^{\frac{5}{2}(1-\sigma)}+r^\sigma K^{1-2\sigma+\epsilon}|t|^{2\sigma}.
\]
\end{proposition}

To prove Proposition \ref{P:1}, a key observation is that the main term arises from the first part of the diagonal term, while the remaining terms are cancelled by a residue from the off-diagonal term.

In the next step, we will deduce an asymptotic formula for the mollified second moment of $L(s,\mathrm{sym}^2f)$ in weight aspect. 
To begin the mollification process, we first need to choose the mollifier $M(s,\mathrm{sym}^2f)$.
Let $P(x)=\frac{B_{2x}(7,7)}{B(7,7)}$, where 
\[
	B_x(k,r)=\int_{0}^{x}u^{k-1}(1-u)^{r-1}du,
\]
and $B(k,r)=B_1(k,r)$ is the Beta function. We observe that $P(x)$ satisfies $P(0)=P'(0)=P'(1/2)=0$ and $P(1/2)=1$.
Next we define the cut-off function. Let $M=K^\theta$ with $0<\theta<1/8$, $s(n)$ is the squarefree kernel of $n$ and
\[
F(x) = 
\begin{cases}
    1 & \text{if } 0\leq x\leq \sqrt{M}, \\
    P\left(\frac{\log(M/x)}{\log M}\right) & \text{if } \sqrt{M}
\leq x\leq M,\\
0& \text{if } x\geq M.
\end{cases}
\]
Then $F$ satiefies $F(m)=0$ if $m>M$. We choose the mollifier \[
M(s,\mathrm{sym}^2f)=\sum_{m\geq1}\frac{b_f(m)F(s(m))}{m^s}.
\]
The mollified second moment takes the following shape:

\begin{proposition}\label{P:2} Let $0<|t|\leq K^{\delta}$ with $0<\delta<1/5$.
When $\sigma>3/2-\epsilon$, we have
\[
	M(\sigma+it,\mathrm{sym}^2f)L(\sigma+it,\mathrm{sym}^2f)=1+O\left(K^{\left(-\frac{1}{20}+\epsilon\right)\theta\left(\sigma-\frac{1}{2}\right)}\right).
\]
When $1/2+1/\log K<\sigma<3/2-\epsilon$, we have
\begin{align*}
&\quad\sum_{k\equiv 0 \bmod 2}\psi\left(\frac{k-1}{K}\right)\sideset{}{^h}\sum_{f\in H_k} |M(\sigma+it,\mathrm{sym}^2f)L(\sigma+it,\mathrm{sym}^2f)|^2\\
&=\sum_{k\equiv 0 \bmod 2}\psi\left(\frac{k-1}{K}\right)+O(K^{1+\left(-\frac{1}{10}+\epsilon\right)(\sigma-\frac{1}{2})}).
\end{align*}
\end{proposition}

\section{Preliminaries}

We need the Petersson trace formula and the approximate functional equation as following: 

\begin{lemma}[Petersson trace formula]\label{L:1} We have
\[
\sideset{}{^h}\sum_{f \in H_k} \lambda_f(n) \lambda_f(m) = \delta_{m,n} + 2\pi i^k \sum_{c=1}^{\infty} \frac{S(n, m; c)}{c} J_{k-1}\left(\frac{4\pi \sqrt{mn}}{c}\right),
\]
where $\delta_{m,n}=1$ if $m=n$ and 0 otherwise, $S(m,n;c)=\sum_{ad\equiv1\bmod c}e\left(\frac{ma+nd}{c}\right)$ is the Kloosterman sum and $J_{k-1}(x)=\sum_{l\geq0}\frac{(-1)^l}{l!\Gamma(l+k)}\left(\frac{x}{2}\right)^{k-1+2l}$ is the $J$-Bessel function.\\
\end{lemma}
\begin{proof} See Proposition 14.5 in \cite{iwaniec2004analytic}.
\end{proof}

\begin{lemma}[Approximate functional equation]\label{L:2} We have
\[
|L(\sigma+it,\mathrm{sym}^2f)|^2=\sum_{m=1}^{\infty}\frac{\lambda_f(m^2)}{m^{\sigma+it}}
\sum_{n=1}^{\infty}\frac{\lambda_f(n^2)}{n^{\sigma-it}}\left (  V_{k,\sigma,t}(mn)+\frac{1}{(mn)^{1-2\sigma}}\tilde{V}_{k,\sigma,t}(mn) \right ),
\]
where, for $\xi>0$ and $G(s)=e^{s^2}$, we have
\begin{align*}
V_{k,\sigma,t}(\xi)&=\frac{1}{2\pi i}\int_{(3)}\frac{L_{\infty}(\sigma+it+s,\mathrm{sym}^2f)L_{\infty}(\sigma-it+s,\mathrm{sym}^2f)}{L_{\infty}(\sigma+it,\mathrm{sym}^2f)L_{\infty}(\sigma-it,\mathrm{sym}^2f)}\\
&\quad\times\zeta(2\sigma+2it+2s)\zeta(2\sigma-2it+2s)\frac{1}{\xi^s}\frac{G(s)}{s}ds,
\end{align*}
and
\begin{align*}
\tilde{V}_{k,\sigma,t}(\xi)&=\frac{1}{2\pi i}\int_{(3)}\frac{L_{\infty}(1-\sigma+it+s,\mathrm{sym}^2f)L_{\infty}(1-\sigma-it+s,\mathrm{sym}^2f)}{L_{\infty}(\sigma+it,\mathrm{sym}^2f)L_{\infty}(\sigma-it,\mathrm{sym}^2f)}\\&\quad\times\zeta(2-2\sigma+2it+2s)\zeta(2-2\sigma-2it+2s)\frac{1}{\xi^s}\frac{G(s)}{s}ds.
\end{align*}
\end{lemma}
\begin{proof}
 Consider the integral
\[
I(\sigma+it,\mathrm{sym}^2f)=\frac{1}{2\pi i}\int_{(3)}\Lambda(\sigma+it+s,\mathrm{sym}^2f)\Lambda(\sigma-it+s,\mathrm{sym}^2f)
G(s)\frac{ds}{s}.
\]
Then we follow Theorem 5.3 in \cite{iwaniec2004analytic} to finish the proof.
\end{proof}

The following Lemma seperates the dependence on $k$ from the gamma functions in $V_{k,\sigma,t}$ and $\tilde{V}_{k,\sigma,t}$.

\begin{lemma}\label{L:3} $V_{k,\sigma,t}$ and $\tilde{V}_{k,\sigma,t}$ are real-valued and satisfiy the following:
\[
V_{k,\sigma,t}(\xi)=\frac{1}{2\pi i}\int_{(3)}\frac{(k^2)^s}{(4\pi^3)^s}\zeta(2\sigma+2it+2s)\zeta(2\sigma-2it+2s)\frac{1}{\xi^s}\frac{G_+(s,\sigma,t)}{s}ds+O(k^{-1+\epsilon}),
\]
where
\[
G_+(s,\sigma,t)=G(s)\frac{\Gamma((\sigma+it+s+1)/2)\Gamma((\sigma-it+s+1)/2)}{\Gamma((\sigma+it+1)/2)\Gamma((\sigma-it+1)/2)}.
\]
\begin{align*}
    \tilde{V}_{k,\sigma,t}(\xi)&=\frac{1}{2\pi i}\int_{(3)}\frac{(k^2)^{1-2\sigma+s}}{(4\pi^3)^{1-2\sigma+s}}\zeta(2-2\sigma+2it+2s)\zeta(2-2\sigma-2it+2s) \\
    & \quad\times\frac{1}{\xi^s}\frac{G_-(s,\sigma,t)}{s}ds +O(k^{-2\sigma+\epsilon}),
\end{align*}
where
\[
G_-(s,\sigma,t)=G(s)\frac{\Gamma((2-\sigma+it+s)/2)\Gamma((2-\sigma-it+s)/2)}{\Gamma((\sigma+it+1)/2)\Gamma((\sigma-it+1)/2)}.
\]
\end{lemma}
\begin{proof} By Stirling's formula, for $|\Delta|<\Re(z)^{1/2}$ , we use the estimate
\[
\frac{\Gamma(z+\Delta)}{\Gamma(z)}=\exp\left(\Delta\log z+\frac{\Delta^2}{2z}+O(|\Delta z|^{-1}) \right).
\]
Hence if $|\Im(s)|<k^\epsilon$, 
\begin{align*}
\frac{\Gamma(\sigma+it+s+k-1)\Gamma(\sigma-it+s+k-1)}{\Gamma(\sigma+it+k-1)\Gamma(\sigma-it+k-1)}&=(k^2+t^2)^s(1+O(k^{-1+\epsilon}))\\
&=k^{2s}(1+O(k^{-1+\epsilon}+|s|t^2/k^2)).
\end{align*}
The contribution of $|\Im(s)|>k^\epsilon$ is negligible since $G(s)$ is a rapid decay function. By Stirling's formula, we can get $G_+(s,\sigma,t)\ll 
    G(s)|t|^s$ and $G_-(s,\sigma,t)\ll 
    G(s)|t|^{1/2-\sigma+s}$. 
Then we complete the proof directly by moving the contour.
\end{proof}

By Petersson trace formula, we need to take an average over $k$ to deal with $J$-Bessel functions.

\begin{lemma}[Average of the $J$-Bessel function]\label{L:4} For $t>0$, we have 
\begin{align*}
\sum_{k \equiv 0 \bmod 2} 2i^{k} h\left(\frac{k-1}{K}\right) J_{k-1}(t) &= -\frac{K}{\sqrt{t}} \Im \left( e^{-2\pi i/8} e^{it} \hbar\left(\frac{K^2}{2t}\right) \right) \\
&\quad+ O\left( \frac{t}{K^4} \int_{-\infty}^{\infty} v^4 |\hat{h}(v)| \, dv \right),
\end{align*}
where \[
\hbar(v)=\int_{0}^{\infty}\frac{h(\sqrt{u})}{\sqrt{2\pi u}}e^{iuv}du,
\] $\hat{h}$ denotes the Fourier transform of $h$.
\end{lemma}
\begin{proof} See Lemma 2.3 in \cite{khan2010non}.\end{proof}

Choosing $h(u)=\psi(u)V_{uK+1,\sigma,t}(n)$ in Lemma \ref{L:4}, we denote $W_{K,\sigma,t}(n,v)$ as $\hbar(v)$ in Lemma \ref{L:4}, i.e.
\[
W_{K,\sigma,t}(n,v)=\int_{0}^{\infty}\frac{\psi(\sqrt{u})V_{\sqrt{u}K+1,\sigma,t}(n)}{\sqrt{2\pi u}}e^{iuv}du.
\]
We also denote $W_{\sigma,t}(n,v)$ as the main term of $W_{K,\sigma,t}(n,v)$ using Lemma \ref{L:3} for $V$ after some modification, i.e.
\[
W_{\sigma,t}\left(n ,v\right)=\frac{1}{2\pi i}\int _{(A)}\frac{1}{(4\pi^3)^s}\zeta(2\sigma+2it+2s)\zeta(2\sigma-2it+2s)\check{\psi}_{s}(v)\frac{1}{n^s}\frac{G_+(s,\sigma,t)}{s}ds,
\]
where
\[
\check{\psi}_s(v)=\int_{0}^{\infty}\frac{\psi(\sqrt{u})}{\sqrt{2\pi u}}u^se^{iuv}du.
\]
Similarly we denote
\[
\tilde{W}_{K,\sigma,t}(n,v)=\int_{0}^{\infty}\frac{\psi(\sqrt{u})\tilde{V}_{\sqrt{u}K+1,\sigma,t}(n)}{\sqrt{2\pi u}}e^{iuv}du,
\]
\begin{align*}
\tilde{W}_{\sigma,t}\left(n ,v\right)&=\frac{1}{2\pi i}\int _{(A)}\frac{1}{(4\pi^3)^s}\zeta(2-2\sigma+2it+2s)\zeta(2-2\sigma-2it+2s)\check{\psi}_{1-2\sigma+s}(v)\\
&\quad\times\frac{1}{n^s}\frac{G_-(s,\sigma,t)}{s}ds.
\end{align*}
By repeated integration by part we have
\[
\frac{\partial ^j}{\partial v^j}{\check{\psi}}_s(v)\ll_{\Re s,B}(1+j+|s|)^B|v|^{-B},
\]
for $j,B\geq0$. By Lemma \ref{L:3},
\[
W_{K,\sigma,t}(n,v)=W_{\sigma,t}\left(\frac{n}{K^2} ,v\right)+O(K^{-1+\epsilon}),
\]
\[
\tilde{W}_{K,\sigma,t}(n,v)=\left(\frac{K^2}{4\pi ^3}\right)^{1-2\sigma}\tilde{W}_{\sigma,t}\left(\frac{n}{K^2} ,v\right)+O(K^{-2\sigma+\epsilon}).
\]
By integration by part, we have 
\[
W_{K,\sigma,t}(n,v)\ll_{A,B} \frac{K^{2A}|t|^A}{n^Av^B}K^{\epsilon} \text{ for } A>0, B\geq0,
\]
\[
\tilde{W}_{K,\sigma,t}(n,v)\ll_{A,B} \frac{K^{2A}|t|^A}{n^Av^B}K^{2-4\sigma+\epsilon} \text{ for } A>\sigma-1/2, B\geq0. 
\]

We need a Mellin transform of $\check{\psi}_{s}(v)$ defined as above.
\begin{lemma}[Mellin transform of $\check{\psi}_{s}(v)$]\label{L:5}
Denote $\tilde{f}$ the Mellin transform of $f$, then
\[
\tilde{\check{\psi}}_s(x)=\int_0^{\infty}\frac{\psi(\sqrt{u})}{\sqrt{2\pi u}}u^{s-x}\Gamma(x)(\cos(x\pi/2)+i\sin(x\pi/2))du.
\]
It has a bound $ \tilde{\check{\psi}}_{s}(x)\ll_{\Re s}(1+|s|^3)|x|^{-2}$.
\end{lemma}
\begin{proof} See Lemma 3.3 in \cite{khan2010non}.\end{proof}

\section{Twisted second moment}

In this section we calculates a second moment of $L(\sigma+it,\mathrm{sym}^2f)$ twisted by $\lambda_f(r^2)$, where $\sigma$ and $t$ are in the range of Proposition \ref{P:1}.

Recall Hecke relation,\[
\lambda_f(n^2)\lambda_f(r^2)=\sum_{d|(n^2,r^2)}\lambda_f(n^2r^2/d^2).
\]
By approximate functional equation and Petersson trace formula,
\[
\begin{aligned}
    &\quad\sum_{k\equiv 0 \bmod 2}\psi\left(\frac{k-1}{K}\right)\sideset{}{^h}\sum_{f\in H_k} |L(\sigma+it,\mathrm{sym}^2f)|^2\lambda_f(r^2)\\
    &=\sum_{k\equiv 0 \bmod 2}\psi\left(\frac{k-1}{K}\right)\sideset{}{^h}\sum_{f\in H_k} \sum_{m=1}^{\infty}\sum_{n=1}^{\infty}\frac{\lambda_f(m^2)\lambda_f(n^2)\lambda_f(r^2)}{m^{\sigma+it}n^{\sigma-it}}\\
&\quad\times\left (  V_{k,\sigma,t}(mn)+\frac{1}{(mn)^{1-2\sigma}}\tilde{V}_{k,\sigma,t}(mn) \right ) \\
    &=\sum_{k\equiv 0 \bmod 2}\psi\left(\frac{k-1}{K}\right)\sum_{n=1}^{\infty}\sum_{d|(n^2,r^2)}\sum_{\substack{m\geq 1\\m=nr/d}}\frac{1}{m^{\sigma+it}}\frac{1}{n^{\sigma-it}}\\
&\quad\times
    \left (  V_{k,\sigma,t}(mn)+\frac{1}{(mn)^{1-2\sigma}}\tilde{V}_{k,\sigma,t}(mn) \right )\\
    &\quad+\sum_{k\equiv 0 \bmod 2}\psi\left(\frac{k-1}{K}\right)\sum_{m,n,c=1}^{\infty}\sum_{d|(n^2,r^2)}
    2\pi i^k\frac{S(m^2,n^2r^2/d^2,c)}{c}J_{k-1}\left ( \frac{4\pi mnr/d}{c} \right )
    \\
    &\quad\times\frac{1}{m^{\sigma+it}}\frac{1}{n^{\sigma-it}}\left (  V_{k,\sigma,t}(mn)+\frac{1}{(mn)^{1-2\sigma}}\tilde{V}_{k,\sigma,t}(mn) \right ).
\end{aligned}
\]
Write $d=d_1 d_2^2$ with $d_1$ squarefree. Note that $d|n^2$ if and only if $d_1d_2|n$, we change the variable $n$ by $nd_1d_2$:
 \[
\begin{aligned}
    &\quad\sum_{k\equiv 0 \bmod 2}\psi\left(\frac{k-1}{K}\right){\sum_{f\in H_k}}^h |L(\sigma+it,\mathrm{sym}^2f)|^2\lambda_f(r^2)\\
    &=\sum_{k\equiv 0 \bmod 2}\psi\left(\frac{k-1}{K}\right)\sum_{n=1}^{\infty}\sum_{d|r^2}\frac{1}{(n^2d_1r)^{\sigma}}\frac{1}{(r/d)^{it}}\\
&\quad\times
    \left (  V_{k,\sigma,t}(n^2d_1r)+\frac{1}{(n^2d_1r)^{1-2\sigma}}\tilde{V}_{k,\sigma,t}(n^2d_1r) \right )\\
    &\quad+\sum_{k\equiv 0 \bmod 2}\psi\left(\frac{k-1}{K}\right)\sum_{m,n,c=1}^{\infty}\sum_{d|r^2}
    2\pi i^k\frac{S(m^2,n^2r^2/d_2^2,c)}{c}J_{k-1}\left ( \frac{4\pi mnr/d_2}{c} \right )
    \\
    &\quad\times\frac{1}{m^{\sigma+it}}\frac{1}{(nd_1d_2)^{\sigma-it}}\left (  V_{k,\sigma,t}(mnd_1d_2)+\frac{1}{(mnd_1d_2)^{1-2\sigma}}\tilde{V}_{k,\sigma,t}(mnd_1d_2) \right ).
\end{aligned}
\]

\subsection{The diagonal term}

Write $(D)=(D_1)+(D_2)$, where $(D_1)$ contains $V$ and $(D_2)$ contains $\tilde{V}$. by Lemma \ref{L:3},
\begin{align*}
    (D_1)  = &\sum_{k\equiv 0 \bmod 2}\psi\left(\frac{k-1}{K}\right) \sum_{d|r^2}\sum_{n=1}^{\infty}\frac{1}{(n^2d_1r)^{\sigma}
    (r/d)^{it}}\frac{1}{2\pi i}\int_{(3)}\frac{(k^2)^s}{(4\pi^3)^s}\zeta(2\sigma+2it+2s)\\
    &\times\zeta(2\sigma-2it+2s)\frac{1}{(n^2d_1r)^s}\frac{G_+(s,\sigma,t)}{s}ds+O(K^{\epsilon})
    \\
    =&\sum_{k\equiv 0 \bmod 2}\psi\left(\frac{k-1}{K}\right) \sum_{d|r^2}\frac{1}{(d_1r)^{\sigma}(r/d)^{it}}\frac{1}{2\pi i}\int_{(\epsilon)}\frac{(k^2)^s}{(4\pi^3)^s}\zeta(2\sigma+2it+2s)\\
    &\times\zeta(2\sigma-2it+2s)\zeta(2\sigma+2s)\frac{1}{(d_1r)^s}\frac{G_+(s,\sigma,t)}{s}ds+O(K^{\epsilon}). 
\end{align*}
We move the contour to $(-1/2-\sigma/2+\epsilon)$ and yield a residue
\begin{align*}
    &\sum_{k\equiv 0 \bmod 2}\psi\left(\frac{k-1}{K}\right) \sum_{d|r^2}\frac{1}{(d_1r)^{\sigma}(r/d)^{it}}\biggl( \zeta(2\sigma+2it)\zeta(2\sigma-2it) \zeta(2\sigma)\biggr.\\
 &+\frac{(k^2)^{1/2-\sigma}}{2(4\pi^3)^{1/2-\sigma}}\frac{G_+(1/2-\sigma,\sigma,t)}{1/2-\sigma}\zeta(1+2it)\zeta(1-2it)\frac{1}{(d_1r)^{1/2-\sigma}}  \\
    & +\frac{(k^2)^{1/2-\sigma-it}}{2(4\pi^3)^{1/2-\sigma-it}}\frac{G_+(1/2-\sigma-it,\sigma,t)}{1/2-\sigma-it}\zeta(1-4it)\zeta(1-2it)\frac{1}{(d_1r)^{1/2-\sigma-it}}\\
    &\biggl.+\frac{(k^2)^{1/2-\sigma+it}}{2(4\pi^3)^{1/2-\sigma+it}}\frac{G_+(1/2-\sigma+it,\sigma,t)}{1/2-\sigma+it}\zeta(1+4it)\zeta(1+2it)\frac{1}{(d_1r)^{1/2-\sigma+it}}\biggr).
\end{align*}
The error term is
$\ll r^{1-\sigma}K^{-\sigma+\epsilon}|t|^{\frac{5}{2}(1-\sigma)}$. Similarly, $(D_2)$ is equal to
\begin{align*}
    &\sum_{k\equiv 0 \bmod 2}\psi\left(\frac{k-1}{K}\right) \sum_{d|r^2}\frac{1}{(d_1r)^{1-\sigma}(r/d)^{it}}\frac{1}{2\pi i}\int_{(\sigma-1/2+\epsilon)}\frac{(k^2)^{1-2\sigma+s}}{(4\pi^3)^{1-2\sigma+s}} \\
    &\times\zeta(2-2\sigma+2it+2s)\zeta(2-2\sigma-2it+2s)\zeta(2-2\sigma+2s)\frac{1}{(d_1r)^s}\frac{G_-(s,\sigma,t)}{s}ds\\
&+O(K^{\epsilon}r^{2\sigma-1}).
\end{align*}
We move the contour to $(\sigma-3/2+\epsilon)$ and yield a residue
\begin{align*}
    &\sum_{k\equiv 0 \bmod 2}\psi\left(\frac{k-1}{K}\right) \sum_{d|r^2}\frac{1}{(d_1r)^{1-\sigma}(r/d)^{it}}\\
    &\times\left( \frac{(k^2)^{1-2\sigma}}{(4\pi^3)^{1-2\sigma}}G_-(0,\sigma,t)\zeta(2-2\sigma+2it)\zeta(2-2\sigma-2it) \zeta(2-2\sigma)\right.\\
    &+\frac{(k^2)^{1/2-\sigma}}{2(4\pi^3)^{1/2-\sigma}}\frac{G_-(\sigma-1/2,\sigma,t)}{\sigma-1/2}\zeta(1+2it)\zeta(1-2it)\frac{1}{(d_1r)^{\sigma-1/2}}  \\
    & +\frac{(k^2)^{1/2-\sigma-it}}{2(4\pi^3)^{1/2-\sigma-it}}\frac{G_-(\sigma-1/2-it,\sigma,t)}{\sigma-1/2-it}\zeta(1-4it)\zeta(1-2it)\frac{1}{(d_1r)^{\sigma-1/2-it}}\\
    &\left.+\frac{(k^2)^{1/2-\sigma+it}}{2(4\pi^3)^{1/2-\sigma+it}}\frac{G_-(\sigma-1/2+it,\sigma,t)}{\sigma-1/2+it}\zeta(1+4it)\zeta(1+2it)\frac{1}{(d_1r)^{\sigma-1/2+it}}\right).
\end{align*}
The error term is $\ll rK^{-1+\epsilon}|t|^{3/2+\sigma}.$ 

We write the main term of $(D_1)+(D_2)$ as 
\[
\sum_\pm S_\pm+\sum_\pm\sum_{\alpha\in\{-it,0,it\}}S_{\alpha,\pm} .\]
Here the definition of $S_\pm$ are given in Proposition \ref{P:1} and
\begin{align*}
    S_{\alpha,\pm}(r)&= \frac{1}{2}\sum_{k\equiv 0 \bmod 2}\psi\left(\frac{k-1}{K}\right)\sum_{d|r^2}\frac{1}{(d_1r)^{1/2+\alpha}(r/d)^{it}}\left(\frac{k^2}{4\pi^3}\right)^{1/2-\sigma+\alpha}\\
     & \quad\times\frac{G_\pm(\pm(1/2-\sigma)+\alpha,\sigma,t)}{\pm(1/2-\sigma)+\alpha} \underset{s=2\alpha}{\mathrm{Res}} (\zeta(1-2it+s)\zeta(1+s)\zeta(1+2it+s)).
\end{align*}
We observe that if $\alpha=0$, since $G(s)$ is an even function,
\begin{align*}
	&\sum_\pm \frac{G_\pm(\pm(1/2-\sigma),\sigma,t)}{\pm(1/2-\sigma)}\\
=&\sum_\pm G(1/2-\sigma)\frac{\Gamma\left(\frac{3/2+it}{2}\right)\Gamma\left(\frac{3/2-it}{2}\right)}{\Gamma\left(\frac{1+\sigma+it}{2}\right)\Gamma\left(\frac{1+\sigma-it}{2}\right)}\left(\frac{1}{1/2-\sigma}+\frac{1}{\sigma-1/2}\right)=0,
\end{align*}
thus $\sum_\pm S_{0,\pm}$=0.

\subsection{Sum of Kloosterman sums}
We follow Lemma 3.2 in \cite{khan2010non} to remove the Kloosterman sums arise from the Petersson trace formula.

\begin{lemma}\label{L:6} For $r<K^{1-2\delta-\epsilon_1}$ with $0<\epsilon_1<1/200$,
\begin{align*} 
&\quad\sum_{m\geq1}S(m^2,b^2,c)e(2mb/c)\cos \left(t\log\frac{bd}{mr}\right)m^{-(1/2+\sigma)}W_{\sigma,t}\left(\frac{mbd}{rK^2},\frac{K^2c}{8\pi mb} \right)
\\
&=\begin{cases}
    \frac{\phi (c)}{\sqrt{c}}\sum_{m\geq1}\cos \left(t\log\frac{bd}{mr}\right)m^{-(1/2+\sigma)}W_{\sigma,t}\left(\frac{mbd}{rK^2},\frac{K^2c}{8\pi mb} \right)+O(K^{-100}) & \text{if } c=\square, \\
    O(K^{-100}) & \text{otherwise. } 
\end{cases}
\end{align*}
\end{lemma}
\begin{proof}We can assume $mb\leq rK^2|t|\leq K^{3-\delta-\epsilon_1}$ and $c\leq K^{1-\delta-\epsilon_1}$, otherwise it becomes $O(K^{-100})$. We have
\begin{align*}
    &\quad\sum_{m\geq1}S(m^2,b^2,c)e(2mb/c)\cos \left(t\log\frac{bd}{mr}\right)m^{-(1/2+\sigma)}W_{\sigma,t}\left(\frac{mbd}{rK^2},\frac{K^2c}{8\pi mb} \right) \\
    & = \sum_{a\bmod c}S(a^2,b^2,c)e(2ab/c)\sum_{m\equiv a\bmod c}\cos \left(t\log\frac{bd}{mr}\right)m^{-(1/2+\sigma)}W_{\sigma,t}\left(\frac{mbd}{rK^2},\frac{K^2c}{8\pi mb} \right).
\end{align*}
We have a derivative estimate of $W_{\sigma,t}$ for $\xi>0$:
\[
\frac{\partial^j}{\partial \xi^j}W_{\sigma,t}\left(\frac{\xi bd}{rK},\frac{Kc}{8\pi \xi b}
\right)\ll_{j,A,B}\xi^{-j-A}\left (\frac{rK|t|}{bd} \right )^A \left(\frac{\xi b}{K^{1-\epsilon}c} \right )^B \sum_{l=1}^{j}\left(\frac{Kc}{b} \right )^lK^{\epsilon}.
\]
Using this bound, we can extend the definition of $W_{\sigma,t}$ to $\xi\leq 0$ by setting \[W_{\sigma,t}\left(\frac{\xi bd}{rK},\frac{Kc}{8\pi \xi b}\right)=0.\] By Poisson summation on arithmetic progression, the $m$-sum is
\[
\frac{K^{1/2-\sigma}}{c}\sum_{m}\int_{-\infty}^{\infty}\cos \left(t\log\frac{bd}{\xi rK}\right)\xi^{-(1/2+\sigma)}W_{\sigma,t}\left(\frac{\xi bd}{rK},\frac{Kc}{8\pi \xi b}
 \right)
  e\left(-\frac{m\xi K}{c} \right)e\left(-\frac{ma}{c} \right)d\xi.
\]
By integration by part $j$ times, if $m\neq 0$, the integral is \[
\ll_{j,A,B}\left ( \frac{c|t|}{mK} \right )^j\left ( \frac{rK|t|}{bd} \right )^A\left ( \frac{b}{Kc} \right )^{B-j}K^{10+B\epsilon}\ll_{\epsilon_1} K^{-100},
\]\\
by choosing $A=1$ and $B=j$ large and $\epsilon_1>2\epsilon$. Then only $m=0$ term survive. By Poisson summation again, the main term equals
\[
\frac{1}{c}\sum_{m}\cos \left(t\log\frac{bd}{mr}\right)m^{-(1/2+\sigma)}W_{\sigma,t}\left(\frac{mbd}{rK^2},\frac{K^2c}{8\pi mb} \right).
\]

Same to Lemma 3.2 in \cite{khan2010non},
\[
\sum_{a\bmod c}S(a^2,b^2,c)e(2ab/c) = 
\begin{cases}
    \phi (c)\sqrt{c} & \text{if } c\  \text{is a square,} \\
    0 & \text{otherwise. } 
\end{cases}
\]
In other cases, we just need to replace $\cos$ by $\sin$ or to replace $W$ by $\tilde{W}$ and make a slight modification to complete the proof.\end{proof}

\subsection{The off-diagonal term}
For the term containing $V$, we let $ \tau=4\pi mnr/cd_2 $,
\begin{align*}
(OD_1)&=\sum_{k\equiv 0 \bmod 2}\psi\left(\frac{k-1}{K}\right)
    \sum_{m,n,c=1}^{\infty}\sum_{d|r^2}
    2\pi i^k\frac{S(m^2,n^2r^2/d_2^2,c)}{c}J_{k-1}(\tau)\\
    &\times\frac{1}{m^{\sigma+it}}\frac{1}{(nd_1d_2)^{\sigma-it}}
    V_{k,\sigma,t}(mnd_1d_2).
\end{align*}
Choosing $
h(\frac{k-1}{K})=\psi(\frac{k-1}{K})V_{k,\sigma,t}(mnd_1d_2)
$ in Lemma \ref{L:4}, we have
\[
\hbar\left( \frac{K^2}{2\tau}\right)=\int_{0}^{\infty}\frac{\psi(\sqrt{u})V_{\sqrt{u}K+1,\sigma,t}(mnd_1d_2)}{\sqrt{2\pi u}}e^{iu\frac{K^2}{2\tau}}du=W_{K,\sigma,t}\left(mnd_1d_2,\frac{K^2cd_2}{8\pi mnr}\right).
\]
Hence by Lemma \ref{L:4}, $(OD_1)$ equals to
\begin{align*}
    &-\pi K\sum_{m,n,c=1}^{\infty}\sum_{d|r^2}\frac{S(m^2,n^2r^2/d_2^2,c)}{cm^{\sigma+it}(nd_1d_2)^{\sigma-it}}\\
&\times
\Im\left( \sqrt{\frac{cd_2}{4\pi mnr}}e\left(\frac{2mnr}{cd_2}-\frac{1}{8} \right)W_{K,\sigma,t}\left(mnd_1d_2,\frac{K^2cd_2}{8\pi mnr}\right)\right) \\
    &+O\left(K^{-4}\sum_{m,n,c=1}^{\infty}\sum_{d|r^2}\frac{r (mn)^{1-\sigma}}{d_1^{\sigma}d_2^{1+\sigma}}\frac{|S(m^2,n^2r^2/d_2^2,c)|}{c^2}\int_{-\infty}^{\infty}v^4|\hat{h}(v)|dv\right).
\end{align*}
By integration by part, the error term is
\[
O\left(K^{-4}\sum_{m,n,c=1}^{\infty}\sum_{d|r^2}\frac{r (mn)^{1-\sigma}}{d_1^{\sigma}d_2^{1+\sigma}}\frac{|S(m^2,n^2r^2/d_2^2,c)|}{c^2}\frac{K^{2A}|t|^A}{(mnd_1d_2)^{A}}K^\epsilon\right),
\]
for any $A>0$. By Weil's bound of Kloosterman sum, $|S(m,n,c)|\ll c^{1/2+\epsilon}(m,n,c)^{1/2}$,
hence the error term is $O(rK^{-2\sigma+\epsilon}|t|^{2-\sigma})$. Similarly, the error term arises from $W_{K,\sigma,t}$
is 
\begin{align*}
& O\left( K^{-1+\epsilon} \sum_{mn\leq K^{2+\epsilon}|t|}\sum_{c\leq nmr/K^{2-\epsilon}}\sum_{d|r^2}\frac{|S(m^2,n^2r^2/d_2^2,c)|}{c(mnd_1d_2)^{\sigma}}\sqrt{\frac{cd_2}{4\pi mnr}}\right)\\
=&O(K^{-2\sigma+\epsilon}|t|^{3/2-\sigma}\sqrt{r}).
\end{align*}

The remaining main term is 
\begin{align*}
    &-\frac{\sqrt{\pi}}{2\sqrt{r}} K\sum_{m,n,c=1}^{\infty}\sum_{d|r^2}\frac{S(m^2,n^2r^2/d_2^2,c)}{c^{1/2}m^{1/2+\sigma+it}n^{1/2+\sigma-it}d_1^{\sigma-it}d_2^{-1/2+\sigma-it}}
\\
&\times\Im\left( e\left(\frac{2mnr}{cd_2}-\frac{1}{8} \right)W_{\sigma,t}\left(\frac{mnd_1d_2}{K^2},\frac{K^2cd_2}{8\pi mnr}\right)\right) \\
     =&-\frac{\sqrt{\pi}}{2\sqrt{r}} K\sum_{m,n,c=1}^{\infty}\sum_{d|r^2}
    \left(\cos \left(t\log\frac{(nd_1d_2)}{m}\right)+i\sin \left(t\log\frac{(nd_1d_2)}{m} \right)\right)\\
   &\times \Im\left(\frac{S(m^2,n^2r^2/d_2^2,c)}{c^{1/2}m^{1/2+\sigma}n^{1/2+\sigma}d_1^{\sigma}d_2^{-1/2+\sigma}} 
 e\left(\frac{2mnr}{cd_2}-\frac{1}{8} \right)W_{\sigma,t}\left(\frac{mnd_1d_2}{K^2},\frac{K^2cd_2}{8\pi mnr}\right)\right).
\end{align*}
Choosing $b=nr/d_2$ in Lemma \ref{L:6}, we have the main term is 
\begin{align*}
&-\frac{\sqrt{\pi}}{2\sqrt{r}} K\Im\biggl(\frac{1}{2}\sum_\pm\sum_{m,n,c=1}^{\infty}\sum_{d|r^2}\frac{\phi(c^2)}{c^2m^{1/2+\sigma\pm it}n^{1/2+\sigma\mp it}d_1^{\sigma\mp it}d_2^{-1/2+\sigma\mp it}}
 e\left(-\frac{1}{8} \right)\biggr.\\
&\times\left.W_{\sigma,t}\left(\frac{mnd_1d_2}{K^2},\frac{K^2c^2d_2}{8\pi mnr}\right)\right) \\
&-i\frac{\sqrt{\pi}}{2\sqrt{r}} K\Im\biggl(\frac{1}{2i}\sum_\pm\sum_{m,n,c=1}^{\infty}\sum_{d|r^2}\frac{\pm\phi(c^2)}{c^2m^{1/2+\sigma\pm it}n^{1/2+\sigma\mp it}d_1^{\sigma\mp it}d_2^{-1/2+\sigma\mp it}}
 e\left(-\frac{1}{8} \right)\biggr.\\
&\times\left.W_{\sigma,t}\left(\frac{mnd_1d_2}{K^2},\frac{K^2c^2d_2}{8\pi mnr}\right)\right).
\end{align*}

We denote the two terms as $(CT_1)$ and $(ST_1)$ respectively. By definition of $W_{\sigma,t}$, we have $(CT_1)$ is equal to
\begin{align*}
&-\frac{\sqrt{\pi}}{2\sqrt{r}}K\Im\left(\frac{1}{2}\sum_\pm\sum_{m,n,c=1}^{\infty}\sum_{d|r^2}\frac{\phi(c^2)}{c^2m^{1/2+\sigma\pm it}n^{1/2+\sigma\mp it}d_1^{\sigma\mp it}d_2^{-1/2+\sigma\mp it}}e(-1/8)\right.\\
&\times\frac{1}{(2\pi i)^2}\int_{(A)}\int_{(B)}\frac{1}{(4\pi^3)^s} \zeta(2\sigma+2it+2s)\zeta(2\sigma-2it+2s)\\
     & \left.\times 
  \left(\frac{K^2c^2d_2}{8\pi mnr} \right)^{-x}\left(\frac{K^2}{mnd_1d_2} \right)^s
  \tilde{\check{\psi}}_s(x)\frac{G_+(s,\sigma,t)}{s}dxds
  \right).
\end{align*}
Choosing $A=1$ and $B=1-\epsilon$ so that $1/2+\sigma+it+s-x$ is in convergence range of zeta function, the above term becomes
\begin{align*}
&-\frac{\sqrt{\pi}}{2\sqrt{r}}K\Im\left(e(-1/8)\frac{1}{(2\pi i)^2}\int_{(1)}\int_{(1-\epsilon)}
(8\pi)^x\frac{G_+(s,\sigma,t)}{s}\zeta(2\sigma+2it+2s)\right.\\
&\times \zeta(2\sigma-2it+2s)\frac{K^{2s-2x}}{(4\pi^3)^s}\zeta(1/2+\sigma+it+s-x)\zeta(1/2+\sigma-it+s-x)\\
&\times\biggl.\zeta(2x)\zeta(2x+1)^{-1}\tilde{\check{\psi}}_s(x)\sum_{d|r^2}\frac{r^x\cos(t\log(d_1d_2))}{d_1^s d_2^{s+x}d_1^{\sigma}d_2^{-1/2+\sigma}}dxds \biggr).
\end{align*}
Moving the contour from $\Re(s)=1$ to $\Re(s)=\epsilon$, crossing simple poles at $s=1/2-\sigma+x\pm it$, we get the residue
\begin{align*}
&\sum_{\pm}-\frac{\sqrt{\pi}}{2\sqrt{r}}K\Im\Biggl(e(-1/8)\frac{1}{2\pi i}\int_{(1-\epsilon)}
(8\pi)^x\frac{G_+(1/2-\sigma+x\pm it,\sigma,t)}{1/2-\sigma+x\pm it}\zeta(1+2x\pm4it)\Biggr.\\
&\times\frac{K^{1-2\sigma\pm2it}}{(4\pi^3)^{1/2-\sigma+x\pm it}}\zeta(1\pm2it)\zeta(2x)\tilde{\check{\psi}}_{1/2-\sigma+x\pm it}(x)\\
&\Biggl.
\times\sum_{d|r^2}\frac{r^x\cos(t\log(d_1d_2))}{d_1^{1/2-\sigma+x\pm it} d_2^{1/2-\sigma+2x\pm it}d_1^{\sigma}d_2^{-1/2+\sigma}}dx \Biggr).
\end{align*}
The error term from $(\epsilon)$ is $O(K^{-1+\epsilon}r^{1/2}|t|^{1-\sigma})$ by moving the $x$-contour to $\Re(x)=1/2+\epsilon$.
Noting that $ e(-1/8)\left(\cos\frac{\pi x}{2}+i\sin\frac{\pi x}{2} \right)=\cos\left(\frac{\pi x}{2}-\frac{\pi}{4}\right)+i\sin\left(\frac{\pi x}{2}-\frac{\pi}{4}\right) $, and by Lemma \ref{L:5}, the above term is equal to
\begin{align*}
&\Re\left(-\frac{\sqrt{\pi}}{\sqrt{r}}\sum_{d|r^2}\frac{1}{2\pi i}\int_{(1-\epsilon)}\frac{K^{2-2\sigma+2it}}{(4\pi^3)^{1/2-\sigma+x+ it}}\sin\left(\frac{\pi x}{2}-\frac{\pi}{4}\right)\right.
(8\pi)^x\Gamma(x)\\
&\times\frac{G_+(1/2-\sigma+x+ it,\sigma,t)}{1/2-\sigma+x+ it}\zeta(1+2x+4it)\zeta(1+2it)\zeta(2x)\int_{0}^{\infty}\frac{\psi(\sqrt{u})}{\sqrt{2\pi u}}u^{
1/2-\sigma+ it}du\\
&\left.\times\frac{r^x\cos(t\log(d_1d_2))}{d_1^{1/2+x+ it} d_2^{2x+ it}}dx \right).
\end{align*}

We can treat $(ST_1)$ similarly, therefore we have
\begin{align*}
(OD_1)=&\sum_{d|r^2}(d_1d_2)^{it}\Re\left(-\frac{\sqrt{\pi}}{\sqrt{r}}\frac{1}{2\pi i}\int_{(1-\epsilon)}\frac{K^{2-2\sigma+2it}}{(4\pi^3)^{1/2-\sigma+x+ it}}\sin\left(\frac{\pi x}{2}-\frac{\pi}{4}\right)\right.
(8\pi)^x\Gamma(x)\\
&\times\frac{G_+(1/2-\sigma+x+ it,\sigma,t)}{1/2-\sigma+x+ it}\zeta(1+2x+4it)\zeta(1+2it)\zeta(2x)\\
&\times\left.\int_{0}^{\infty}\frac{\psi(\sqrt{u})}{\sqrt{2\pi u}}u^{1/2-\sigma+ it}du
\frac{r^x}{d_1^{1/2+x+ it} d_2^{2x+ it}}dx \right).
\end{align*}

Similarly, we consider the term containing $\tilde{V}$ and get 
\begin{align*}
(OD_2)=& \sum_{d|r^2}(d_1d_2)^{it}\Re\left(-\frac{\sqrt{\pi}}{\sqrt{r}}\frac{1}{2\pi i}\int_{(1-\epsilon)}\frac{K^{2-2\sigma+2it}}{(4\pi^3)^{1/2-\sigma+x+ it}}\sin\left(\frac{\pi x}{2}-\frac{\pi}{4}\right)
(8\pi)^x\Gamma(x)\right.\\
&\times\frac{G_-(\sigma-1/2+x+ it,\sigma,t)}{\sigma-1/2+x+ it}\zeta(1+2x+4it)\zeta(1+2it)\zeta(2x)
 \\
&\left.\times\int_{0}^{\infty}\frac{\psi(\sqrt{u})}{\sqrt{2\pi u}}u^{
1/2-\sigma+ it}du\frac{r^x}{d_1^{1/2+x+it} d_2^{2x+it}}dx\right)\\
&+O(rK^{-1+\epsilon}|t|^{1+\sigma}+r^{\sigma}K^{1-2\sigma+\epsilon}|t|^{2\sigma}).
\end{align*}
We denote the main term of $(OD_1)+(OD_2)$ as $\sum_\pm I_\pm$, here \begin{align*}
I_\pm(r)=&\sum_{d|r^2}(d_1d_2)^{it}\Re \left(-\frac{\sqrt{\pi}}{\sqrt{r}}\frac{1}{2\pi i}\int_{(1-\epsilon)}\frac{K^{2-2\sigma+2it}}{(4\pi^3)^{1/2-\sigma+x+ it}}\sin\left(\frac{\pi x}{2}-\frac{\pi}{4}\right)(8\pi)^x\Gamma(x)\right.\\
&\times\frac{G_\pm(\pm(1/2-\sigma)+x+ it,\sigma,t)}{\pm(1/2-\sigma)+x+ it}\zeta(1+2x+4it)\zeta(1+2it)\zeta(2x)\\
&\left.\int_{0}^{\infty}\frac{\psi(\sqrt{u})}{\sqrt{2\pi u}}u^{1/2-\sigma+ it}du\frac{r^x}{d_1^{1/2+x+it} d_2^{2x+it}}dx\right).
\end{align*}
\subsection{Cancellation}
\begin{lemma}\label{L:7} Assume $f(K,t)$ is a positive, monotonically decreasing function of $K$, satisfying $ \lim_{K \to \infty}f(K,t)=0 $. $g(t)$ is a function valued on $\mathbb{C}$. If for any sufficiently large $K$ and $0<t<K^\delta$, $ \Im(K^{it}g(t))\ll f(K,t) $, then $g(t)=0$.
\end{lemma}
\begin{proof} Let $t$ be fixed, we can find a sequence $\{K_n\}_{n=n_0}^\infty$ where \[K_n=\exp\left(\frac{1}{t}\left(\frac{\pi}{2}+n\pi-\arg g(t)\right)\right).\] Therefore we have that there exists a constant $c$, $\Im(K_n^{it}g(t))=|K_n^{it}g(t)|=|g(t)|\leq cf(K_n,t)$. Then we let $n\rightarrow\infty$ to complete the proof.
\end{proof}

By the functional equation of zeta function we have $\Gamma(x)\zeta(2x)=\pi^{x-1/2}\Gamma(\frac{1-2x}{2})\zeta(1-2x)$. Applying the duplication and reflection formula of $\Gamma(x)$, we get  
\[
\sin\left(\frac{\pi x}{2}-\frac{\pi}{4}\right)\Gamma(x)\zeta(2x)=-2^{-1/2-x}\pi^{2x}\frac{\Gamma(\frac{3}{4}-\frac{x}{2})}{\Gamma(\frac{3}{4}+\frac{x}{2})}\zeta(1-2x).
\]
Thus $(OD_1)+(OD_2)$ is equal to
\begin{align*}
\sum_\pm I_\pm(r)=&\sum_{d|r^2}(d_1d_2)^{it}\Re \left(\frac{1}{2\sqrt{r}}\frac{1}{2\pi i}\int_{(1-\epsilon)}\frac{K^{2-2\sigma+2it}}{(4\pi^3)^{1/2-\sigma+ it}}\frac{\Gamma(\frac{3}{4}-\frac{x}{2})}{\Gamma(\frac{3}{4}+\frac{x}{2})}\right.\\
&\times\sum_\pm\frac{G_\pm(\pm(1/2-\sigma)+x+ it,\sigma,t)}{\pm(1/2-\sigma)+x+ it}\zeta(1+2x+4it)\zeta(1+2it)\zeta(1-2x)\\
&\times\left.\int_{0}^{\infty}\frac{\psi(\sqrt{u})}{\sqrt{u}}u^{1/2-\sigma+ it}du\frac{r^x}{d_1^{1/2+x+it} d_2^{2x+it}}dx\right).
\end{align*}

An observation is that if we calculate the twisted second moment again, but combining $\lambda_f(m^2)$ and $\lambda_f(r^2)$instead of $\lambda_f(n^2)$ and $\lambda_f(r^2)$ by Hecke relation, we get the same diagonal term since $ \sum_{d|r^2}(r/d)^xd_1^v $ is an even function of $x$, but a different off-diagonal term $(OD_1')+(OD_2')$:
\begin{align*}
\sum_\pm\tilde{I}_\pm(r)=&\sum_{d|r^2}(d_1d_2)^{-it}\Re \left(\frac{1}{2\sqrt{r}}\frac{1}{2\pi i}\int_{(1-\epsilon)}\frac{K^{2-2\sigma+2it}}{(4\pi^3)^{1/2-\sigma+ it}}\frac{\Gamma(\frac{3}{4}-\frac{x}{2})}{\Gamma(\frac{3}{4}+\frac{x}{2})}\right.\\
&\times\sum_\pm\frac{G_\pm(\pm(1/2-\sigma)+x+ it,\sigma,t)}{\pm(1/2-\sigma)+x+ it}\zeta(1+2x+4it)\zeta(1+2it)\zeta(1-2x)\\
&\times\left.\int_{0}^{\infty}\frac{\psi(\sqrt{u})}{\sqrt{u}}u^{1/2-\sigma+ it}du\frac{r^x}{d_1^{1/2+x+it} d_2^{2x+it}}dx\right).
\end{align*}
Then we have the difference of $(OD_1)+(OD_2)$ and $(OD_1')+(OD_2')$ can be bounded by the error term arise from the off-diagonal term, since we are calculating the same second moment. 

Writing \[\sum_\pm\tilde{I}_\pm(r)=\sum_{d|r^2}(d_1d_2)^{-it}\Re(K^{2-2\sigma+2it}J(t,d_1,d_2)),\] we have \[\sum_\pm I_\pm(r)-\sum_\pm\tilde{I}_\pm(r)=O(R')=2i\sum_{d|r^2}\Re\left(K^{2-2\sigma+2it}J(t,d_1,d_2)\sin(t\log d_1d_2)\right),\]
where $R'= K^{-1+\epsilon}|t|^{1+\sigma}r+K^{1-2\sigma+\epsilon}|t|^{2\sigma}r^\sigma$. When $\sigma<3/2-\epsilon$ for some $\epsilon>0$, we choose $f(K,t)=R'K^{2\sigma-2}$ in Lemma \ref{L:7}, then we have\[\sum_{d|r^2}J(t,d_1,d_2)(d_1d_2)^{it}=\sum_{d|r^2}J(t,d_1,d_2)(d_1d_2)^{-it}.\]
Replacing the variable $x$ into $x-it$ and changing the variable $x$ into $-x$ for $I_-$, we have $\sum_\pm I_\pm(r)$ is equal to
\begin{align*}
&\Re \left(\frac{1}{\sqrt{r}}\frac{1}{2\pi i}\int_{(1-\epsilon)}\frac{K^{2-2\sigma+2it}}{(4\pi^3)^{1/2-\sigma+ it}}\frac{\Gamma\left(\frac{3/2+it+x}{2}\right)\Gamma\left(\frac{3/2+it-x}{2}\right)}{\Gamma\left(\frac{\sigma+it+1}{2}\right)\Gamma\left(\frac{\sigma-it+1}{2}\right)}\frac{G(1/2-\sigma+x)}{1/2-\sigma+x}\right.\\
&\times\left.\tilde{\psi}(-2\sigma+2it)\zeta(1+2x+2it)\zeta(1+2it)\zeta(1-2x+2it)\sum_{d|r^2}\frac{(r/d)^{x}(d_1d_2)^{-it}}{r^{it}d_1^{1/2} d_2^{-it}}dx\right)\\
&+\Re\left(-\frac{1}{\sqrt{r}}\frac{1}{2\pi i}\int_{(-1+\epsilon)}\frac{K^{2-2\sigma+2it}}{(4\pi^3)^{1/2-\sigma+ it}}\frac{\Gamma\left(\frac{3/2+it+x}{2}\right)\Gamma\left(\frac{3/2+it-x}{2}\right)}{\Gamma\left(\frac{\sigma+it+1}{2}\right)\Gamma\left(\frac{\sigma-it+1}{2}\right)}\frac{G(-(1/2-\sigma)-x)}{-(1/2-\sigma)-x}\right.\\
&\times\left.\tilde{\psi}(-2\sigma+2it)\zeta(1+2x+2it)\zeta(1+2it)\zeta(1-2x+2it)\sum_{d|r^2}\frac{(r/d)^{-x}(d_1d_2)^{-it}}{r^{it}d_1^{1/2} d_2^{-it}}dx\right),
\end{align*}
where $ \tilde{\psi}(-2\sigma+2it)=\int_{0}^{\infty}\psi(u)u^{1-2\sigma+2it}du $.
Again, we note that $\sum_{d|r^2}(r/d)^xd_1^v $ is an even function of $x$, and $G(s)$ is also an even function of $s$, hence the integrand in $I_+$ and $I_-$ is same. By Cauchy residue theorem we encounter poles at $x=\pm it$ and $x=\sigma-1/2$. At $x=\pm it$ we yield residues 
\begin{align*}
&\sum_\pm\mp\frac{1}{2}\Re \frac{K^{2-2\sigma+2it}}{(4\pi^3)^{1/2-\sigma+ it}}\frac{\Gamma\left(\frac{3/2+2it}{2}\right)\Gamma\left(3/4\right)}{\Gamma\left(\frac{\sigma+it+1}{2}\right)\Gamma\left(\frac{\sigma-it+1}{2}\right)}\frac{G(1/2-\sigma\pm it)}{1/2-\sigma\pm it}\\
&\times\tilde{\psi}(-2\sigma+2it)\zeta(1+4it)\zeta(1+2it)\sum_{d|r^2}\frac{(r/d)^{it}}{r^{1/2+it}d_1^{1/2+it} }.
\end{align*}

Then we turn to $S_{\alpha,\pm}$. By Poisson summation we have \begin{align*}
&\sum_{k\equiv 0 \bmod 2}\psi\left(\frac{k-1}{K}\right)k^{1-2\sigma+2it}\\
=&\frac{K^{2-2\sigma+2it}}{2}\sum_m(-1)^m\int_{0}^{\infty}\psi(u)u^{1-2\sigma+2it}e\left(-\frac{mKu}{2}\right)du+O(K^{1-2\sigma}).
\end{align*}
By stationary phase we have
 \[
\sum_{k\equiv 0 \bmod 2}\psi\left(\frac{k-1}{K}\right)k^{1-2\sigma+2it}=\frac{K^{2-2\sigma+2it}}{2}\int_{0}^{\infty}\psi(u)u^{1-2\sigma+2it}du+O(K^{1-2\sigma}).
\]
Then we get $S_{\alpha,\pm}$ satisfies \begin{align*}
    &\sum_\pm\sum_{\alpha\in\{-it,it\}}S_{\alpha,\pm}\\
 = &\Re\sum_{k\equiv 0 \bmod 2}\psi\left(\frac{k-1}{K}\right)\sum_{d|r^2}\frac{1}{(d_1r)^{1/2+it}(r/d)^{it}}\left(\frac{k^2}{4\pi^3}\right)^{1/2-\sigma+it}\\
     & \times \sum_\pm\frac{G_\pm(\pm(1/2-\sigma)+it,\sigma,t)}{\pm(1/2-\sigma)+it}\zeta(1+2it)\zeta(1+4it)\\
     =&\frac{1}{2}\Re \frac{K^{2-2\sigma+2it}}{(4\pi^3)^{1/2-\sigma+it}}\tilde{\psi}(-2\sigma+2it)\frac{\Gamma\left(\frac{3/2+2it}{2}\right)\Gamma\left(3/4\right)}{\Gamma\left(\frac{\sigma+it+1}{2}\right)\Gamma\left(\frac{\sigma-it+1}{2}\right)}\sum_\pm \frac{G(1/2-\sigma\pm it)}{\pm(1/2-\sigma)+it}\\
     &\times\zeta(1+2it)\zeta(1+4it)\sum_{d|r^2}\frac{1}{(d_1r)^{1/2+it}(r/d)^{it}}+O(K^{1-2\sigma+\epsilon}),
\end{align*}
this term cancels out with the residue at $\pm it$. For the pole at $\sigma-1/2$ we yield a residue 
\begin{align*}
S_{OD}=&\Re\frac{K^{2-2\sigma+2it}}{(4\pi^3)^{1/2-\sigma+ it}}\frac{\Gamma\left(\frac{2-\sigma+it}{2}\right)}{\Gamma\left(\frac{\sigma-it+1}{2}\right)}\tilde{\psi}(-2\sigma+2it)\sum_{d|r^2}\frac{(r/d)^{\sigma-1/2}}{r^{1/2+it}d_1^{1/2+it}}\\
&\times\zeta(2\sigma+2it)\zeta(1+2it)\zeta(2-2\sigma+2it).
\end{align*}
Then Proposition \ref{P:1} follows. \qed
\section{Mollification}

Recall the definition of $F$ defined in introduction, we see that $F$ is continuously differentiable and $F$ has a Mellin transform
\[
\hat{F}(s)=\frac{M^s}{s} {_{1}F_1\left(7;14;-\frac{s\log M}{2}\right)},
\]
where \[
_{1}F_1(k;k+r;b)=\frac{1}{B(k,r)}\int_{0}^{1}u^{k-1}(1-u)^{r-1}e^{bu}du,
\]
is the Kummer's function (see (3.1.3) in \cite{slater1960confluent} for example). By Python we have 
$_{1}F_1(7;14;b)=\frac{P(-b)-e^{-b}P(b)}{b^{13}}$,
where $P(b)$ is given as following:
\begin{align*}
&P(b)=8648640b^6+363242880b^5+7264857600b^4\\
&+87178291200b^3+653837184000b^2+2876883609600b+5753767219200.
\end{align*}
Then we get \[
	\hat{F}(s)=O\left(\frac{M^s+M^{s/2}}{s(s\log M)^7}+...+\frac{M^s+M^{s/2}}{s(s\log M)^{13}}\right).
\]

Let $n$ be a 4th-power-free integer, then $n$ admits a unique decomposition
$n=n_1n_2^2n_3^3$, where $n_1,n_2,n_3$ are squarefree and pairwise coprime. For such $n$ we have
\[
b_f(n)=\mu(n_1)\lambda_f(n_1^2)\lambda_f(n_2^2)\mu(n_3),
\]
and we may therefore rewrite $M(s,\mathrm{sym}^2f)$ as
\[
M(s,\mathrm{sym}^2f)=\sideset{}{^{\flat}}\sum_{n_1,n_2,n_3}\frac{\mu(n_1)\lambda_f(n_1^2)\lambda_f(n_2^2)\mu(n_3)}{n_1^s n_2^{2s} n_3^{3s}}F(n_1n_2n_3).
\]
Here we denote $ {\sum}^\flat $ for the summation over squarefree numbers, and $\sum_{a_1,\cdots,a_n}^*$ for the summation over pairwise coprime $a_1,\cdots,a_n$.
Expanding the square and applying the Hecke relation, we have
\begin{align*}
   &|M(\sigma+it,\mathrm{sym}^2f)|^2\\
=&\sideset{}{^{*\flat}}\sum_{m_1,m_2,m_3}\ \ \sideset{}{^{*\flat}}\sum_{n_1,n_2,n_3}\frac{\mu(m_1)\lambda_f(m_1^2)\lambda_f(m_2^2)\mu(m_3)\mu(n_1)\lambda_f(n_1^2)\lambda_f(n_2^2)\mu(n_3)}{m_1^{\sigma+it} m_2^{2\sigma+2it} m_3^{3\sigma+3it} n_1^{\sigma-it} n_2^{2\sigma-2it} n_3^{3\sigma-3it}}  \\
   &\times F(m_1m_2m_3)F(n_1n_2n_3)\\
   =&\sideset{}{^{*\flat}}\sum_{m_1,m_2,m_3}\ \ \sideset{}{^{*\flat}}\sum_{n_1,n_2,n_3}\frac{\mu(m_1)\mu(m_3)\mu(n_1)\mu(n_3)F(m_1m_2m_3)F(n_1n_2n_3)}{m_1^{\sigma+it} m_2^{2\sigma+2it} m_3^{3\sigma+3it} n_1^{\sigma-it} n_2^{2\sigma-2it} n_3^{3\sigma-3it}} \\
   &\times \sum_{g|(m_1^2m_2^2,n_1^2n_2^2)}\lambda_f\left(\frac{m_1^2m_2^2n_1^2n_2^2}{g^2}\right).
\end{align*}
Since $(m_1,m_2)=(n_1,n_2)=1$, we set $\alpha=(m_1m_2,n_1n_2)=\alpha_{11}\alpha_{12}\alpha_{21}\alpha_{22}$ which satisfies $\alpha_{ij}=(m_i,n_j)$ for $i,j=1,2$. 

By Proposition \ref{P:1} we have
\begin{align*}
   & \sum_{k\equiv 0 \bmod 2}\psi\left(\frac{k-1}{K}\right)\sideset{}{^h}\sum_{f\in H_k}|LM(\sigma+it,\mathrm{sym}^2f)|^2
\\=&\sideset{}{^{*\flat}}\sum_{\substack{\alpha_{11},\alpha_{12},\alpha_{21},\alpha_{22}\\m_1,m_2,n_1,n_2}}\sideset{}{^{\flat}}\sum_{\substack{m_3,n_3\\(m_3n_3,\alpha_{j_1,j_2})=1\\(m_3,m_1m_2)=1\\(n_3,n_1n_2)=1}}
\frac{\mu(m_1n_1)\mu(m_3n_3)\mu(\alpha_{12}\alpha_{21})}{m_1^{\sigma+it} m_2^{2\sigma+2it} m_3^{3\sigma+3it} n_1^{\sigma-it} n_2^{2\sigma-2it} n_3^{3\sigma-3it}} \\
    & \times \sum_{g|\alpha_{11}^2\alpha_{12}^2\alpha_{21}^2\alpha_{22}^2}\frac{F(\alpha_{11}\alpha_{12}\alpha_{21}\alpha_{22}m_1m_2m_3)F(\alpha_{11}\alpha_{12}\alpha_{21}\alpha_{22}n_1n_2n_3)}{\alpha_{11}^{2\sigma}\alpha_{12}^{3\sigma-it}\alpha_{21}^{3\sigma+it}\alpha_{22}^{4\sigma}}
    \\
&\times\sum_{k\equiv 0 \bmod 2}\psi\left(\frac{k-1}{K}\right)\sideset{}{^h}\sum_{f\in H_k} |L(\sigma+it,\mathrm{sym}^2f)|^2\lambda_f\left(\frac{\alpha_{11}^4\alpha_{12}^4\alpha_{21}^4\alpha_{22}^4m_1^2n_1^2m_2^2n_2^2}{g^2}\right)\\
=&\mathscr{S}_++\mathscr{S}_-+\mathscr{S}_{OD}+\mathscr{R},
\end{align*}
where $\mathscr{S}_+$, $\mathscr{S}_-$, $\mathscr{S}_{OD}$, $\mathscr{R}$ are corresponding to $S_+$, $S_-$, $S_{OD}$, $R$ in Proposition \ref{P:1} respectively. 
We let $\tau_{u,v}(r)=\sum_{d|r^2}1/(d_1^u d^v)$ if $d$ has a squarefree decomposition $d=d_1d_2^2$, we can see $\tau_{u,v}$ is a multiplicative function of $n$.

 For the error term $\mathscr{R}$, we have \[
\mathscr{R}\ll \sum_{m_1m_2m_3\leq M}\sum_{n_1n_2n_3\leq M}\sum_{g|(m_1m_2,n_1n_2)^2}\frac{R(m_1n_1m_2n_2/g)}{m_1^\sigma m_2^{2\sigma}m_3^{3\sigma}n_1^\sigma n_2^{2\sigma}n_3^{3\sigma}}.
\]
We can prove\[
	\sum_{m_1m_2m_3\leq M}\frac{1}{m_1^\sigma m_2^{2\sigma}m_3^{3\sigma}}=O(M^{\max\{2-2\sigma,0\}}).
\]
By the bounds $R(m_1n_1m_2n_2/g)\ll M^{4\sigma-2}K^\epsilon$ and $\tau((m_1m_2,n_1n_2)^2)\ll K^\epsilon$, we have\[
	\mathscr{R}\ll M^2K^\epsilon \ll K^{1-\frac{\theta}{2}(\sigma-1/2)},
\]
if $\theta\leq1/5-\epsilon$ and $\delta\leq1/5$.

Before beginning, we need some results for the ratios of $\zeta(s)$.

\begin{lemma}\label{L:8}
(1) For $1/\log K<\delta<1/2$ or $-c_1/(\log K)^{3/4}<|\delta|<-c_2/\log K$, and $t\ll K^a$ we have  \[
\frac{\zeta(1+\delta+it)}{\zeta(1+2\delta+it)}=O(1).
\]  
(2) Assume $1/\log K<\delta<1$ or $-c_1/(\log K)^{3/4}<\delta<-c_2/\log K$. For $|t|<K^{1/5}$ and $l\geq 7$, \[
\int_{-\infty}^{\infty}\left|\frac{\zeta(1+\delta+it)}{\zeta(1+\delta+it+is)}\right|^2\frac{ds}{|\zeta(1+\delta+is)|^2|\delta+is|^{2l}}=O((\log K)^{2l-3}).
\]
\end{lemma}
\begin{proof} For (1), if $\delta>0$ we take a logarithm:\[
\log\frac{\zeta(1+2\delta+it)}{\zeta(1+\delta+it)}=\int_{1+\delta}^{1+2\delta}\frac{\zeta'}{\zeta}(x+it)dx\ll\delta\left|\frac{\zeta'}{\zeta}(1+\delta)\right|\ll1.
\]
If $\delta<0$ we apply the upper bound $\frac{\zeta'}{\zeta}(\sigma+it)\ll (\log t)^{2/3+\epsilon}$ (see section 6.19 in \cite{titchmarsh1986theory}) for $\sigma+it$ in the Vinogradov zero-free region.

For (2), first we write the integrand as $g(\delta,t,s)$. We suppose $\delta<1/(\log K)^{2/3}$ (otherwise is trivial), and divide it into 2 cases:\\
(i) If $|t|\geq\frac{1}{(\log K)^{2/3}}$, when $|s|<\frac{1}{2(\log K)^{2/3}}$, we take a logarithm:
\[
\log\frac{\zeta(1+\delta+it)}{\zeta(1+\delta+it+is)}=-\int_{0}^{s}\frac{\zeta'}{\zeta}(1+\delta+it+ix)dx\ll\delta\left|\frac{\zeta'}{\zeta}(1+it)\right|\ll1.
\]
Therefore we have $\zeta(1+\delta+it)\asymp \zeta(1+\delta+it+is)$, then \[
	\int_0^{1/(\log K)^{2/3}}g(\delta,t,s)ds\ll \int_0^\delta \delta^{2-2l}ds+\left|\int_\delta^{1/(\log K)^{2/3}}|s|^{2-2l}ds\right|\ll (\log K)^{2l-3}.
\]
When $\frac{1}{2(\log K)^{2/3}}<|s|<1$, we have if $l\geq 7$,\begin{align*}
	\int_{1/(\log K)^{2/3}}^1g(\delta,t,s)ds\ll \int_{1/(\log K)^{2/3}}^1\frac{(\log K)^{10/3+\epsilon}}{|\delta+is|^{2l-2}}ds\ll& (\log K)^{\frac{10}{3}+\epsilon+\frac{2}{3}(2l-3)}\\
\ll&(\log K)^{2l-3}.
\end{align*}
When $|s|>1$, we have if $l\geq5$,
 \[
	\int_1^\infty g(\delta,t,s)ds \ll \int_1^\infty(\log K)^6\frac{ds}{|s|^{2l}}\ll (\log K)^{2l-3}.
\]
(ii) If $|t|<\frac{1}{(\log K)^{2/3}}$, when $|s|<1$, we have
\[
	\int_0^1 g(\delta,t,s)ds\ll \int_0^1\frac{|\delta+it+is|^2}{|\delta+it|^2}\frac{ds}{|\delta+is|^{2l-2}}\ll (\log K)^{2l-3},
\]
by discussing the size of $|t|$, $|s|$ and $\delta$.
When $|s|\geq 1$, we can use the same estimate as above. 
\end{proof}
\subsection {The treatment of $\mathscr{S}_+$}
We will show that $\mathscr{S}_+$ contributes the main term of the mollified second moment in this subsection.
By Mellin inversion, we have $\mathscr{S}_+$ is equal to
\begin{align*}
&\sum_{k\equiv 0 \bmod 2}\psi\left(\frac{k-1}{K}\right)\sideset{}{^{*\flat}}\sum_{\substack{\alpha_{11},\alpha_{12},\alpha_{21},\alpha_{22}\\m_1,m_2,n_1,n_2}}\sideset{}{^{\flat}}\sum_{\substack{m_3,n_3\\(m_3n_3,\alpha_{j_1,j_2})=1\\(m_3,m_1m_2)=1\\(n_3,n_1n_2)=1}}
\frac{\mu(m_1n_1)\mu(m_3n_3)\mu(\alpha_{12}\alpha_{21})}{m_1^{\sigma+it} m_2^{2\sigma+2it} m_3^{3\sigma+3it} n_1^{\sigma-it} n_2^{2\sigma-2it} n_3^{3\sigma-3it}} \\
    & \times \frac{1}{(2\pi i)^2}\int_{(c)}\int_{(c)}\frac{\hat{F}(\beta)\hat{F}(\gamma)\alpha_{11}^{-\beta-\gamma}\alpha_{12}^{-\beta-\gamma}\alpha_{21}^{-\beta-\gamma}\alpha_{22}^{-\beta-\gamma}m_1^{-\beta}m_2^{-\beta}m_3^{-\beta}n_1^{-\gamma}n_2^{-\gamma}n_3^{-\gamma}}{\alpha_{11}^{2\sigma}\alpha_{12}^{3\sigma-it}\alpha_{21}^{3\sigma+it}\alpha_{22}^{4\sigma}(\alpha_{11}^2\alpha_{12}^2\alpha_{21}^2\alpha_{22}^2m_1m_2n_1n_2)^{\sigma+it}}d\beta d\gamma
    \\&\times\tau_{\sigma,-it}(m_1m_2n_1n_2)\sum_{g_{11}|\alpha_{11}^2}\sum_{g_{12}|\alpha_{12}^2}\sum_{g_{21}|\alpha_{21}^2}\sum_{g_{22}|\alpha_{22}^2}\tau_{\sigma,-it}\left(\frac{\alpha_{11}^2}{g_{11}}\right)\tau_{\sigma,-it}\left(\frac{\alpha_{12}^2}{g_{12}}\right)\\ 
&\times\tau_{\sigma,-it}\left(\frac{\alpha_{21}^2}{g_{21}}\right)\tau_{\sigma,-it}\left(\frac{\alpha_{22}^2}{g_{22}}\right)(g_{11}g_{12}g_{21}g_{22})^{\sigma+it}\zeta(2\sigma+2it)\zeta(2\sigma-2it)\zeta(2\sigma).
\end{align*}
Denote
\begin{align*}
    &\mathscr{G}(\beta,\gamma,\sigma,t)\\
=&\prod_p \left(1-\frac{1}{p^{2\sigma+2it+\beta}}-\frac{1}{p^{2\sigma-2it+\gamma}}-\frac{1}{p^{2\sigma+\beta}}-\frac{1}{p^{2\sigma+\gamma}}+\frac{1}{p^{2\sigma+\beta+\gamma}}+\frac{1}{p^{4\sigma+2it+\beta}} \right.\\
     & \left.+\frac{1}{p^{4\sigma-2it+\gamma}}+\frac{1}{p^{4\sigma+\beta+\gamma}}-\frac{1}{p^{5\sigma-it+\beta+\gamma}}-\frac{1}{p^{5\sigma-3it+\beta+\gamma}}-\frac{1}{p^{5\sigma+it+\beta+\gamma}}\right.\\
     &\left.-\frac{1}{p^{5\sigma+3it+\beta+\gamma}}+\frac{1}{p^{5\sigma+3it+\gamma}}+\frac{1}{p^{5\sigma+it+\gamma}}+\frac{1}{p^{5\sigma-it+\beta}}+\frac{1}{p^{5\sigma-3it+\beta}}-\frac{1}{p^{6\sigma+\gamma}}-\frac{1}{p^{6\sigma+\beta}}\right.\\
&\left.+\frac{1}{p^{6\sigma+\beta+\gamma}}-\frac{1}{p^{7\sigma+it+\gamma}}-\frac{1}{p^{7\sigma-it+\beta}}+\frac{1}{p^{7\sigma+it+\beta+\gamma}}+\frac{1}{p^{7\sigma-it+\beta+\gamma}}\right).
\end{align*}
Then we can rewrite the inner sum as an Euler product:
\begin{align*}
\mathscr{S}_+=&\zeta(2\sigma+2it)\zeta(2\sigma-2it)\zeta(2\sigma)\sum_{k\equiv 0 \bmod 2}\psi\left(\frac{k-1}{K}\right)\\
&\times\frac{1}{(2\pi i)^2}
\int_{(c)} \int_{(c)}\mathscr{G}(\beta,\gamma,\sigma,t) \hat{F}(\beta)\hat{F}(\gamma)d\beta d\gamma.
\end{align*}
Let the region $\mathscr{D}$ be
\[
\begin{cases}
    \Re(\beta+2\sigma)>1/2,\\
    \Re(\gamma+2\sigma)>1/2, \\
    \Re(\beta+\gamma+2\sigma)>1/2.
\end{cases}
\]
We have \[
\mathscr{G}(\beta,\gamma,\sigma,t)=\frac{\zeta(2\sigma+\beta+\gamma)}{\zeta(2\sigma+\beta)\zeta(2\sigma+\gamma)\zeta(2\sigma+2it+\beta)\zeta(2\sigma-2it+\gamma)}H(\beta,\gamma,\sigma,t),
\]
where $H$ converges absolutely in the region $\mathscr{D}$ and satisfies $H(0,0,\sigma,t)=1$. We move the both integral to $\Re(\beta)=\Re(\gamma)=1/\log K$ and truncate the both integral at $|\Im(\beta)|<K,|\Im(\gamma)|<K$ with an error $O(K^{-1+\epsilon})$.
Recall the zero-free region of the Riemann zeta function, we move the $\beta$-contour to 
\begin{align*}
&\mathscr{C}=\biggl\{\beta:\Re(\beta)\biggr.\\
=&\biggl.1-2\sigma-\frac{c_0}{\log^{3/4}\left((|\Im(\beta)+2t|+10)(|\Im(\beta)-2t|+10)(|\Im(\beta)|+10)\right)}\biggr\},
\end{align*}
and we encounter poles at $\beta=0$ and $\beta=1-2\sigma-\gamma$. The remaining integral is $ O(K^{1-\theta(\sigma-1/2)}e^{-c_1(\log K)^{1/5}}) $.

For the first pole, we yield a residue  \[
\sum_{k\equiv 0 \bmod 2}\psi\left(\frac{k-1}{K}\right)\frac{1}{2\pi i}\int_{1/\log K-iK}^{1/\log K+iK}\frac{\zeta(2\sigma-2it)}{\zeta(2\sigma-2it+\gamma)}H(0,\gamma,\sigma,t)\hat{F}(\gamma)d\gamma.
\]
Moving the contour to $\mathscr{C}\cap\{\gamma:|\Im(\gamma)|<K\}$, we encounter a pole at $\gamma=0$, hence this term contributes\[
\sum_{k\equiv 0 \bmod 2}\psi\left(\frac{k-1}{K}\right)+O(K^{1-\frac{\theta}{2}(\sigma-1/2)}e^{-c_1(\log K)^{1/5}}).
\]

For the second pole, we yield a residue
\begin{align*}
&\sum_{k\equiv 0 \bmod 2}\psi\left(\frac{k-1}{K}\right)\frac{1}{2\pi i}\int_{1/\log K-iK}^{1/\log K+iK}\frac{\zeta(2\sigma+2it)\zeta(2\sigma-2it)\zeta(2\sigma)}{\zeta(1-\gamma)\zeta(2\sigma+\gamma)\zeta(1+2it-\gamma)\zeta(2\sigma-2it+\gamma)}\\
&\times H(1-2\sigma-\gamma,\gamma,\sigma,t)\hat{F}(\gamma)\hat{F}(1-2\sigma-\gamma)d\gamma.
\end{align*}
Moving the contour to $\Re(\gamma)=-\sigma+1/2$, we do not encounter any pole since the pole of $\hat{F}(\gamma)$ cancels with zero of $1/\zeta(1-\gamma)$, then the integral is bounded by
\begin{align*}
&\sum_{k\equiv 0 \bmod 2}\psi\left(\frac{k-1}{K}\right)\int_{-K}^{K}\zeta(2\sigma)\left|\frac{\zeta(2\sigma+2it)\hat{F}(1/2-\sigma+is)}{\zeta(1/2+\sigma+is)\zeta(1/2+\sigma-2it+is)}\right|^2ds\\
\ll&\sup_{8\leq l\leq 14}\sum_{k\equiv 0 \bmod 2}\psi\left(\frac{k-1}{K}\right)\log K|\zeta(2\sigma+2it)|^2 \\
&\times\int_{-K}^{K}\left|\frac{M^{-(\sigma-1/2)/2}|\sigma-1/2+is|^{-l}}{\zeta(1/2+\sigma+is)\zeta(1/2+\sigma-2it+is)(\log M)^{l-1}}\right|^2ds,
\end{align*}
in spite of an error term $O(K^{-1+\epsilon})$. This contributes $O(K^{1-\frac{\theta}{2}(\sigma-\frac{1}{2})})$ by Lemma \ref{L:8}. 

\subsection{The treatment of $\mathscr{S}_-$}
We will show that $\mathscr{S}_-$ contributes an error term. 
We have $\mathscr{S}_-$ is equal to
\begin{align*}
&\sum_{k\equiv 0 \bmod 2}\psi\left(\frac{k-1}{K}\right)\sideset{}{^{*\flat}}\sum_{\substack{\alpha_{11},\alpha_{12},\alpha_{21},\alpha_{22}\\m_1,m_2,n_1,n_2}}\sideset{}{^{\flat}}\sum_{\substack{m_3,n_3\\(m_3n_3,\alpha_{j_1,j_2})=1\\(m_3,m_1m_2)=1\\(n_3,n_1n_2)=1}}
\frac{\mu(m_1n_1)\mu(m_3n_3)\mu(\alpha_{12}\alpha_{21})}{m_1^{\sigma+it} m_2^{2\sigma+2it} m_3^{3\sigma+3it} n_1^{\sigma-it} n_2^{2\sigma-2it} n_3^{3\sigma-3it}} \\
    & \times \frac{F(\alpha_{11}\alpha_{12}\alpha_{21}\alpha_{22}m_1m_2m_3)F(\alpha_{11}\alpha_{12}\alpha_{21}\alpha_{22}n_1n_2n_3)}{\alpha_{11}^{2\sigma}\alpha_{12}^{3\sigma-it}\alpha_{21}^{3\sigma+it}\alpha_{22}^{4\sigma}}\left(\frac{k^2}{4\pi^3}\right)^{1-2\sigma}|\zeta(2-2\sigma+2it)|^2
    \\&\times\zeta(2-2\sigma)G_-(0,\sigma,t)\sum_{g|\alpha_{11}^2\alpha_{12}^2\alpha_{21}^2\alpha_{22}^2}\frac{\tau_{1-\sigma,-it}\left(\alpha_{11}^2\alpha_{12}^2\alpha_{21}^2\alpha_{22}^2m_1m_2n_1n_2/g\right)}{(\alpha_{11}^2\alpha_{12}^2\alpha_{21}^2\alpha_{22}^2m_1m_2n_1n_2/g)^{1-\sigma+it}}.
\end{align*}
We only consider $\sigma-1/2<\frac{200\log\log K}{\log K}$, otherwise $\log K$ can be bounded by $K^{\frac{1}{200}(\sigma-1/2)}$. Then we have a trivival estimate by the convexity bound of the zeta function:
\[
\mathscr{S}_-\ll K^{3-4\sigma} (\log K)^{100} |\zeta(2-2\sigma+2it)|^2 G_-(0,\sigma,t)\ll K^{1-\theta(\sigma-1/2)}.
\] 
Let the region $\mathscr{D}$ be
\[
\begin{cases}
    \Re(\beta)>-1/4,\\
    \Re(\gamma)>-1/4, \\
    \Re(\beta+\gamma)>-1/4.
\end{cases}
\]
We can take a Mellin transform and rewrite the inner sum as an Euler product:
\begin{align*}
&\mathscr{G}(\beta,\gamma,\sigma,t)\\
=& 
 \prod_p\left(1+\frac{1}{p^{2\sigma+\beta+\gamma}}-\frac{1}{p^{1+2it+\beta}}-\frac{1}{p^{1+\beta}}-\frac{1}{p^{1+\gamma}}-\frac{1}{p^{1-2it+\gamma}}+O\left(\frac{1}{p^{2-\sigma+\min\{\beta,\gamma\}}}\right)\right) \\
=&\frac{\zeta(2\sigma+\beta+\gamma)H(\beta,\gamma,\sigma,t)}{\zeta(1+2it+\beta)\zeta(1+\beta)\zeta(1+\gamma)\zeta(1-2it+\gamma)},
\end{align*}
where $H$ converges absolutely in the region $\mathscr{D}$. Noting that $\zeta(2-2\sigma)\ll\log K$,  we move the both integral to $\Re(\beta)=\Re(\gamma)=1/\log K$ and let
 \[
\mathscr{C}=\left\{\beta:\Re(\beta)=-\frac{c_0}{\log^{3/4}\left((|\Im(\beta)+2t|+10)(|\Im(\beta)|+10)(|\Im(\beta)-2t|+10)\right)}\right\}.
\]
Moving the $\beta$-contour to $\mathscr{C}$, we observe that the remaining integral is small, and we encounter a pole at $\beta=1-2\sigma-\gamma$. 

For the pole, we yield a residue \begin{align*}
	&\sum_{k\equiv 0 \bmod 2}\psi\left(\frac{k-1}{K}\right)\frac{1}{2\pi i}\int_{(c)}\frac{\hat{F}(\gamma)\hat{F}(1-2\sigma-\gamma)\zeta(2-2\sigma)|\zeta(2-2\sigma+2it)|^2d\gamma}{\zeta(2-2\sigma-\gamma)\zeta(2-2\sigma+2it-\gamma)\zeta(1-2it+\gamma)\zeta(1+\gamma)}\\
&\times\left(\frac{k^2}{4\pi^3}\right)^{1-2\sigma}G_-(0,\sigma,t).
\end{align*}
Moving the contour to $\Re(\gamma)=-\sigma+1/2$ and truncating at $|\Im\gamma|\leq K$, we do not encounter poles, and the integral contributes $O(K^{1-\frac{\theta}{2}(\sigma-\frac{1}{2})})$ by Lemma \ref{L:8}.
\subsection{The treatment of $\mathscr{S}_{OD}$}
We have $\mathscr{S}_{OD}$ equals to
\begin{align*}
&\Re\frac{K^{2-2\sigma+2it}}{(4\pi^3)^{1/2-\sigma+it}}\sideset{}{^{*\flat}}\sum_{\substack{\alpha_{11},\alpha_{12},\alpha_{21},\alpha_{22}\\m_1,m_2,n_1,n_2}}\sideset{}{^{\flat}}\sum_{\substack{m_3,n_3\\(m_3n_3,\alpha_{j_1,j_2})=1\\(m_3,m_1m_2)=1\\(n_3,n_1n_2)=1}}
\frac{\mu(m_1n_1)\mu(m_3n_3)\mu(\alpha_{12}\alpha_{21})}{m_1^{\sigma+it} m_2^{2\sigma+2it} m_3^{3\sigma+3it} n_1^{\sigma-it} n_2^{2\sigma-2it} n_3^{3\sigma-3it}} \\
    & \times \frac{1}{(2\pi i)^2}\int_{(c)}\int_{(c)}\frac{\hat{F}(\beta)\hat{F}(\gamma)\alpha_{11}^{-\beta-\gamma}\alpha_{12}^{-\beta-\gamma}\alpha_{21}^{-\beta-\gamma}\alpha_{22}^{-\beta-\gamma}m_1^{-\beta}m_2^{-\beta}m_3^{-\beta}n_1^{-\gamma}n_2^{-\gamma}n_3^{-\gamma}}{\alpha_{11}^{2\sigma}\alpha_{12}^{3\sigma-it}\alpha_{21}^{3\sigma+it}\alpha_{22}^{4\sigma}(\alpha_{11}^2\alpha_{12}^2\alpha_{21}^2\alpha_{22}^2m_1m_2n_1n_2)^{1-\sigma+it}}d\beta d\gamma
    \\&\times\tau_{1/2+it,\sigma-1/2}(m_1m_2n_1n_2)\sum_{g_{11}|\alpha_{11}^2}\sum_{g_{12}|\alpha_{12}^2}\sum_{g_{21}|\alpha_{21}^2}\sum_{g_{22}|\alpha_{22}^2}\tau_{1/2+it,\sigma-1/2}\left(\frac{\alpha_{11}^2}{g_{11}}\right)\\
    &\times\tau_{1/2+it,\sigma-1/2}\left(\frac{\alpha_{12}^2}{g_{12}}\right)\tau_{1/2+it,\sigma-1/2}\left(\frac{\alpha_{21}^2}{g_{21}}\right)\tau_{1/2+it,\sigma-1/2}\left(\frac{\alpha_{22}^2}{g_{22}}\right)\frac{\Gamma\left(\frac{2-\sigma+it}{2}\right)}{\Gamma\left(\frac{1+\sigma-it}{2}\right)}\\ 
&\times\tilde{\psi}(-2\sigma+2it)(g_{11}g_{12}g_{21}g_{22})^{1-\sigma+it}\zeta(2\sigma+2it)\zeta(1+2it)\zeta(2-2\sigma+2it),
\end{align*}
where we can also restrict $\sigma$ in the range $\sigma-1/2<\frac{200\log\log K}{\log K}$. We can write the inner sum as an Euler product:
\begin{align*}
&\mathscr{G}(\beta,\gamma,\sigma,t)\\
=& 
 \prod_p\left(1+\frac{1}{p^{2\sigma+\beta+\gamma}}-\frac{1}{p^{1+2it+\beta}}-\frac{1}{p^{2\sigma+2it+\beta}}-\frac{1}{p^{1+\gamma}}-\frac{1}{p^{2\sigma+\gamma}}+O\left(\frac{1}{p^{1+\sigma+\min\{\beta,\gamma\}}}\right)\right) \\
=&\frac{\zeta(2\sigma+\beta+\gamma)H(\beta,\gamma,\sigma,t)}{\zeta(1+2it+\beta)\zeta(2\sigma+2it+\beta)\zeta(1+\gamma)\zeta(2\sigma+\gamma)}.
\end{align*}
Let the region $\mathscr{D}$ be
\[
\begin{cases}
    \Re(\beta)>-\sigma+1/4,\\
    \Re(\gamma)>-\sigma+1/4, \\
    \Re(\beta+\gamma)>-\sigma+1/4,
\end{cases}
\]
then $H$ converges absolutely in the region $\mathscr{D}$. We get \begin{align*}
\mathscr{S}_{OD}
=&\Re\frac{K^{2-2\sigma+2it}}{(4\pi^3)^{1/2-\sigma+it}}\frac{\Gamma\left(\frac{2-\sigma+it}{2}\right)}{\Gamma\left(\frac{1+\sigma-it}{2}\right)}\frac{1}{(2\pi i)^2}\int_{(c)} \int_{(c)}\mathscr{G}(\beta,\gamma,\sigma,t) \hat{F}(\beta)\hat{F}(\gamma)d\beta d\gamma\\
&\times \tilde{\psi}(-2\sigma+2it)\zeta(2\sigma+2it)\zeta(1+2it)\zeta(2-2\sigma+2it).
\end{align*}
Again, we move the $\beta$-contour to $\mathscr{C}$ with the same choice as subsection 4.2. We encounter poles at $\beta=0$ and $\beta=1-2\sigma-\gamma$. 

For the first pole we yield a residue
\begin{align*}
\Re\frac{K^{2-2\sigma+2it}}{(4\pi^3)^{1/2-\sigma+it}}\frac{\Gamma\left(\frac{2-\sigma+it}{2}\right)}{\Gamma\left(\frac{1+\sigma-it}{2}\right)}\tilde{\psi}(-2\sigma+2it)\frac{1}{(2\pi i)}
\int_{(c)}\frac{\zeta(2-2\sigma+2it)}{\zeta(1+\gamma)}\hat{F}(\gamma)d\gamma.
\end{align*}
Moving the $\gamma$-contour to $\mathscr{C}$, the pole of $\hat{F}$ at $\gamma=0$ does not contribute because it cancels with the zero of $1/\zeta(1+\gamma)$.

For the residue of the second pole, we truncate the $\gamma$-integral at $|\Im\gamma|\leq K$ and move the contour to $\Re\gamma=-\sigma+1/2$. In spite of an negligible error term, we get the residue is bounded by
\begin{align*}
&\sup_{8\leq l\leq 14}\Re\left(\frac{K^{2-2\sigma+2it}}{(4\pi^3)^{1/2-\sigma+it}}\frac{\Gamma\left(\frac{2-\sigma+it}{2}\right)}{\Gamma\left(\frac{1+\sigma-it}{2}\right)}\tilde{\psi}(-2\sigma+2it)\frac{1}{2\pi i}\int_{-K}^{K}\frac{M^{(1/2-\sigma)/2}}{|1/2-\sigma+is|^{2l}(\log M)^{2l-2}}\right.\\
&\times\left.\left|\frac{\zeta(2\sigma+2it)\zeta(1+2it)\zeta(2-2\sigma+2it)}{\zeta(3/2-\sigma+2it-is)\zeta(1/2+\sigma+2it-is)\zeta(3/2-\sigma+is)\zeta(1/2+\sigma+is)}\right|ds\right).
\end{align*}
We apply Cauchy-Schwarz inequality to obtain that the inner integral is bounded by
	\begin{align*}
	&\left(\int_{-K}^{K}\frac{M^{(1/2-\sigma)/2}|\zeta(1+2it)|}{|1/2-\sigma+is|^{2l}(\log M)^{2l-2}}\left|\frac{\zeta(2\sigma+2it)}{\zeta(1/2+\sigma+2it-is)\zeta(1/2+\sigma+is)}\right|^2ds\right)^{1/2}\\
	&\times\left(\int_{-K}^{K}\frac{M^{(1/2-\sigma)/2}|\zeta(1+2it)|}{|1/2-\sigma+is|^{2l}(\log M)^{2l-2}}\left|\frac{\zeta(2-2\sigma+2it)}{\zeta(3/2-\sigma+2it-is)\zeta(3/2-\sigma+is)}\right|ds\right)^{1/2}.
	\end{align*}

When $t>\frac{1}{100\log K}$, we have $\zeta(1+2it)=O(\log K)$, then $\mathscr{S}_\sigma$ is $O(K^{1-\frac{\theta}{2}(\sigma-\frac{1}{2})})$ by Lemma \ref{L:8}.

When $0<t<\frac{1}{100\log K}$, we follow the method by Conrey-Soundararajan in \cite{conrey2002real}. Let $\mathscr{O}$ be the circle with center 0 and radius $\frac{1}{5\log K}$, then by Cauchy's residue theorem, 
\begin{align*}
	\mathscr{S}_{OD}=&\frac{1}{2\pi i}\int_{\mathscr{O}}\frac{K^{2-2\sigma+2w}}{(4\pi^3)^{1/2-\sigma+w}}\frac{\Gamma\left(\frac{2-\sigma+w}{2}\right)}{\Gamma\left(\frac{1+\sigma-w}{2}\right)}\tilde{\psi}(-2\sigma+2w)\eta(\sigma,w)\\
&\times\zeta(2\sigma+2w)\zeta(1+2w)\zeta(2-2\sigma+2w)\left(\frac{1}{w-it}+\frac{1}{w+it}\right)dw.
\end{align*}
Where
\begin{align*}
&\eta(\sigma,w)\\
=&\sum_{\substack{\alpha_{11},\alpha_{12},\alpha_{21},\alpha_{22}\\m_1,m_2,n_1,n_2}}\sideset{}{^{\flat}}\sum_{\substack{m_3,n_3\\(m_3n_3,\alpha_{j_1,j_2})=1\\(m_3,m_1m_2)=1\\(n_3,n_1n_2)=1}}
\frac{\mu(m_1n_1)\mu(m_3n_3)\mu(\alpha_{12}\alpha_{21})}{m_1^{\sigma+w} m_2^{2\sigma+2w} m_3^{3\sigma+3w} n_1^{\sigma-w} n_2^{2\sigma-2w} n_3^{3\sigma-3w}} \\
    & \times \frac{1}{(2\pi i)^2}\int_{(c)}\int_{(c)}\frac{\hat{F}(\beta)\hat{F}(\gamma)\alpha_{11}^{-\beta-\gamma}\alpha_{12}^{-\beta-\gamma}\alpha_{21}^{-\beta-\gamma}\alpha_{22}^{-\beta-\gamma}m_1^{-\beta}m_2^{-\beta}m_3^{-\beta}n_1^{-\gamma}n_2^{-\gamma}n_3^{-\gamma}}{\alpha_{11}^{2\sigma}\alpha_{12}^{3\sigma-w}\alpha_{21}^{3\sigma+w}\alpha_{22}^{4\sigma}(\alpha_{11}^2\alpha_{12}^2\alpha_{21}^2\alpha_{22}^2m_1m_2n_1n_2)^{1-\sigma+w}}d\beta d\gamma
    \\&\times\tau_{1/2+w,\sigma-1/2}(m_1m_2n_1n_2)\sum_{g_{11}|\alpha_{11}^2}\sum_{g_{12}|\alpha_{12}^2}\sum_{g_{21}|\alpha_{21}^2}\sum_{g_{22}|\alpha_{22}^2}\tau_{1/2+w,\sigma-1/2}\left(\frac{\alpha_{11}^2}{g_{11}}\right)\\
    &\times\tau_{1/2+w,\sigma-1/2}\left(\frac{\alpha_{12}^2}{g_{12}}\right)\tau_{1/2+w,\sigma-1/2}\left(\frac{\alpha_{21}^2}{g_{21}}\right)\tau_{1/2+w,\sigma-1/2}\left(\frac{\alpha_{22}^2}{g_{22}}\right)(g_{11}g_{12}g_{21}g_{22})^{1-\sigma+w}\\
=&\frac{1}{(2\pi i)^2}\int_{(c)}\int_{(c)}\frac{\zeta(2\sigma+\beta+\gamma)\hat{F}(\beta)\hat{F}(\gamma)H(\beta,\gamma,\sigma,-iw)}{\zeta(1+2w+\beta)\zeta(2\sigma+2w+\beta)\zeta(1+\gamma)\zeta(2\sigma+\gamma)}d\beta d\gamma. 
\end{align*}
By the same evaluation to previous cases we have $\mathscr{S}_\sigma=O(K^{1-\frac{\theta}{2}(\sigma-\frac{1}{2})})$, then we complete the proof of Proposition \ref{P:2}.\qed

\section{Proof of Theorem \ref{T:1}}
We follow Selberg's original proof in \cite{selberg1946contributions} to establish the theorem. We first need a Lemma due to Selberg:

\begin{lemma}\label{L:11} Let $f(s)$ be regular in and upon the boundary of the domain $\sigma\geq\alpha$, $t_1\leq t\leq t_2$, where $s=\sigma+it$, further suppose that in this domain
\[
	f(s)=1+o\left(\exp\left(-\frac{\pi}{t_2-t_1}\sigma\right)\right),
\]
uniformly as $\sigma\rightarrow\infty$. Denote zeros of $f(s)$ as $\beta+i\gamma$, then we have
\begin{align*}
&(t_2-t_1) \sum_{\substack{\beta>\sigma\\t_1<\gamma<t_2}}\sin\left(\pi\frac{\gamma-t_1}{t_2-t_1}\right)\left(\exp\left(\pi\frac{\beta-\alpha}{t_2-t_1}\right)-\exp\left(-\pi\frac{\beta-\alpha}{t_2-t_1}\right)\right)\\
=&\int_{t_1}^{t_2}\sin\left(\pi\frac{t-t_1}{t_2-t_1}\right)\log|f(\alpha+it)|dt\\
&+\int_{\sigma'}^{\infty}\left(\exp\left(\pi\frac{\sigma_1-\alpha}{t_2-t_1}\right)-\exp\left(-\pi\frac{\sigma_1-\alpha}{t_2-t_1}\right)\right)\\
&\times\left(\log|f(\sigma+it_1)|+\log|f(\sigma+it_2)|\right)d\sigma_1.
\end{align*}
\end{lemma}
\begin{proof} See Lemma 14 in \cite{selberg1946contributions}.\end{proof}

It suffices to consider the case $t_2-t_1=\frac{1}{\log K}$. Let $\xi=K^{\theta/3}$ and $t_1'=t_1-\frac{2}{\log\xi}$, $t_2'=t_2+\frac{2}{\log\xi}$, $\sigma'=\sigma-\frac{1}{2\log K}$.
Using the inequality $ e^x-e^{-x}>2x $ for $x>0$, we have
\begin{align*}
&N(\sigma,t_1,t_2,\mathrm{sym}^2 f)=\sum_{\substack{\beta>\sigma\\t_1<\gamma<t_2}}1\leq \sum_{\substack{\beta>\sigma'\\t_1'<\gamma<t_2'}}2(\beta-\sigma')\log K\\
&\ll \log K(t_2'-t_1') \sum_{\substack{\beta>\sigma'\\t_1'<\gamma<t_2'}}\sin\left(\pi\frac{\gamma-t_1'}{t_2'-t_1'}\right)\left(\exp\left(\pi\frac{\beta-\sigma'}{t_2'-t_1'}\right)-\exp\left(-\pi\frac{\beta-\sigma'}{t_2'-t_1'}\right)\right).
\end{align*}
Note that if $\rho=\beta+i\gamma$ is a zero of $L(s,\mathrm{sym}^2 f)$, then $\rho$ is a zero of $ML(s,\mathrm{sym}^2 f)$. By Lemma \ref{L:11} and the inequality $\log |x|\leq \frac{|x|^2-1}{2}$, we have
\begin{align*}
N(\sigma,t_1,t_2,\mathrm{sym}^2 f)\ll& \log K \int_{t_1'}^{t_2'}\sin\left(\pi\frac{\gamma-t_1'}{t_2'-t_1'}\right)
\left(|ML(\sigma'+it,\mathrm{sym}^2 f)|-1\right)dt\\
&+\log K\int_{\sigma'}^{\infty}\left(\exp\left(\pi\frac{\sigma_1-\sigma'}{t_2'-t_1'}\right)-\exp\left(-\pi\frac{\sigma_1-\sigma'}{t_2'-t_1'}\right)\right)\\
&\times\left(|ML(\sigma_1+it_1')|^2+|ML(\sigma_1+it_2')|^2-2\right)d\sigma_1.
\end{align*}
Noting that $\exp(\frac{\pi}{t_2'-t_1'})\ll \xi^{\pi/4}$, by Proposition \ref{P:2} we obtain
\begin{align*}
&\sum_{k\equiv0\bmod2}\psi\left(\frac{k-1}{K}\right)\sideset{}{^h}\sum_{f\in H_k}N(\sigma,t_1,t_2,\mathrm{sym}^2 f)\\
\ll& K^{1-\frac{1}{21}(\sigma-\frac{1}{2})}\log K(t_2'-t_1')+\int_{\sigma'}^{\infty}\xi^{\frac{\pi}{4}(\sigma_1-\sigma')}K^{1-\frac{1}{21}(\sigma_1-\frac{1}{2})}\log Kd\sigma_1\ll K^{1-\frac{1}{21}(\sigma-\frac{1}{2})},
\end{align*}
which completes the proof of Theorem \ref{T:1}.\qed

\subsection*{Acknowledgements} I would like to thank Professor Bingrong Huang for helpful discussions.
\subsection* {Funding} This work was supported by the Key R\&D Program of Shandong Province, China (No. 2026CXPT038).
\subsection* {AI Declaration} AI tools were used in choosing a suitable polynomial $P(x)$ in the mollifier $M(s,\mathrm{sym}^2f)$.

\end{document}